\documentclass[12pt,letterpaper]{iopart}

 \usepackage{iopamsm}
 \usepackage[noadjust]{cite}
 \usepackage{color}
 \usepackage{curves}
 \usepackage{enumerate}
 \usepackage{epic}
 \usepackage{epsfig}
 \usepackage{float}
 \usepackage{graphics}
 \usepackage{graphicx}
 \usepackage{latexsym}
 \usepackage{multicol}
 \usepackage{subfigure}
 \usepackage{theorem}

\usepackage{layout}

 \newtheorem{thm}{Theorem}[section]
 \newtheorem{hyp}{Hypothesis}
 
 \newtheorem{lemma}[thm]{Lemma}
 \newtheorem{prop}[thm]{Proposition}
 \theorembodyfont{\upshape}
 \newtheorem{rem}[thm]{Remark}
\numberwithin{equation}{section}
 \DeclareMathOperator{\A}{A}
 
 \DeclareMathOperator{\BB}{B}

 \DeclareMathOperator{\rank}{rank}
 \DeclareMathOperator{\real}{Re}

 \DeclareMathOperator{\T}{T}
 
 \DeclareMathOperator{\GL}{GL}
 
 \newcommand{\proof}{\textbf{Proof:\quad}}

 \newcommand{\GC}{\mathbb{C}}

 \newcommand{\GM}{\mathbb{M}}

 \newcommand{\GR}{\mathbb{R}}

 \newcommand{\GSO}{\mathbb{SO}}

 \newcommand{\ov}{\overline}

 \newcommand{\newl}{\newline\newline}
  \DeclareRobustCommand{\qed}{%
  \ifmmode \mathqed
  \else
    \leavevmode\unskip\penalty9999 \hbox{}\nobreak\hfill
    \quad\hbox{$\Box$\normalsize}%
  \fi
}
\renewcommand{\footnoterule}{%
\kern -3pt 
\hrule height 0.6pt width 0.4\columnwidth
\kern 2.6pt 
} \textheight=630pt
\begin{document}
\title[Spiral anchoring in media with multiple inhomogeneities]{Spiral anchoring in media with multiple inhomogeneities: a dynamical
system approach}

\author{P Boily\footnote{Present address: Institute of the Environment, University of Ottawa,
Ottawa K1N 6N5, Canada.}, V G LeBlanc and E Matsui}

\address{Department of Mathematics and Statistics, University of Ottawa,
Ottawa K1N 6N5, Canada} \ead{pboily@uottawa.ca}
\begin{abstract}
The spiral is one of Nature's more ubiquitous shape: it can be seen in various
media, from galactic geometry to cardiac tissue. In the literature, very
specific models are used to explain some of the observed incarnations of these
dynamic entities. Barkley \cite{B1,B2} first noticed that the range of possible
spiral behaviour is caused by the Euclidean symmetry that these models
possess.\par In experiments however, the physical domain is never perfectly
Euclidean. The heart, for instance, is finite, anisotropic and littered with
inhomogeneities. To capture this loss of symmetry (and as a result model the
physical situation with a higher degree of accuracy), LeBlanc and Wulff
introduced forced Euclidean symmetry-breaking (FESB) in the analysis,
\textit{via} two basic types of perturbations: translational symmetry-breaking
(TSB) and rotational sym\-metry-breaking terms. In \cite{LW,LeB}, they show
that phenomena such as anchoring and quasi-periodic meandering can be explained
by combining Barkley's insight with FESB.
\par In this article, we provide a fuller characterization of spiral anchoring by studying the effects of $n$ simultaneous TSB
perturbations, where $n>1$.
\end{abstract}

\submitto{Journal of Nonlinear Science}
\maketitle


\section{Introduction} \normalsize Spiral waves have been observed
in a variety of experimental contexts, ranging from the well-known
Belousov-Zhabotinsky chemical reaction to the electrical potential in cardiac
tissue \cite{LOPS,GZM,ZM,MPMPV,WR,YP,Detal,R1,Wetal,J,MAK,B1,B2,BKT}. In this
last case, spiral waves are believed to be a precursor to several fatal cardiac
arrythmias ({\it e.g.} ventricular tachycardia and ventricular fibrillation)
\cite{W,Wetal,KS}. A thorough understanding of the various dynamical properties
of spiral waves is therefore warranted.

One of the most interesting and fruitful approaches in recent years to the
study of spiral waves has been to use the theory of equivariant dynamical
systems to derive finite-dimensional models for many of the observed dynamical
states and bifurcations of spirals. The pioneer of this approach was Barkley,
who realized that the experimentally-observed transition from rigid rotation to
quasi-periodic meandering and drifting could be explained using only the
underlying symmetries (the group $\mathbb{S}\mathbb{E}(2)$ of all planar
translations and rotations) of the governing reaction-diffusion partial
differential equations: he derived an {\it ad hoc} system of 5 ordinary
differential equations with $\mathbb{S}\mathbb{E}(2)$ symmetry which model a
Hopf bifurcation from a rotating wave, and then showed that this
finite-dimensional system replicated the experimentally-observed transition to
meandering and drifting \cite{B1,B2,BK}. Sandstede, Scheel and Wulff later
proved a general center manifold reduction theorem for relative equilibria and
relative periodic solutions in spatially extended infinite-dimensional
$\mathbb{S}\mathbb{E}(2)$-equivariant dynamical systems, thereby providing
mathematical justification for Barkley's approach
\cite{SSW1,SSW2,SSW3,SSW4,FSSW}.

One of the advantages of this equivariant dynamical systems approach is that
one can often give universal, model-independent explanations of many of the
observed dynamics and bifurcations of spiral waves. For example, the
above-mentioned Hopf bifurcation from rigid rotation to quasi-periodic
meandering and drifting has been observed in both numerical simulations
\cite{BKT} and in actual chemical reactions \cite{LOPS}. Another example is the
anchoring/repelling of spiral waves on/from a site of inhomogeneity, which has
been observed in numerical integrations of an Oregonator system \cite{MPMPV},
in photo-sensitive chemical reactions \cite{ZM} and in cardiac tissue
\cite{Detal}. Using a model-independent approach based on forced
symmetry-breaking, LeBlanc and Wulff showed that anchoring/repelling of
rotating waves is a generic property of systems in which the translation
symmetry of $\mathbb{S}\mathbb{E}(2)$ is broken by a small perturbation
\cite{LW}. Similarly, some dynamics of spiral waves observed in anisotropic
media ({\em e.g.} phase-locking and/or linear drifting of meandering spiral
waves) have been shown to be generic consequences of rotational
symmetry-breaking \cite{R1,R2,LeB,Bo3,Bo2,Bo1}.

Consider as a paradigm a system of reaction-diffusion partial differential
equations
\begin{equation}
\frac{\partial u}{\partial t}=D\cdot\nabla^2\,u+f(u) \label{basicrdpde}
\end{equation}
where $u$ is a $k$-vector valued function of time and two-dimensional space,
$D$ is a matrix of diffusion coefficients and $f:\mathbb{R}^k\longrightarrow
\mathbb{R}^k$ is a smooth reaction term.  Many of the phenomena in which spiral
waves are observed experimentally are modeled by systems of the form
(\ref{basicrdpde}).  Moreover, Scheel has proved that systems of this form can
admit time-periodic, rigidly rotating spiral wave solutions \cite{S}. Implicit
in the form of equations (\ref{basicrdpde}) is the fact that the medium of
propagation is completely homogeneous and isotropic. Mathematically, this is
represented by the invariance of (\ref{basicrdpde}) under the transformations
\begin{equation}
u(t,x)\longmapsto u(t,x_1\cos\,\theta-x_2\sin\,\theta+p_1,x_1\sin\,\theta+
x_2\cos\,\theta+p_2), \label{uaction}
\end{equation}
where $(\theta,p_1,p_2)\in\mathbb{S}^1\times\mathbb{R}^2$ and $x\in \GR^2$
\cite{Wulff,DMcK}. The group of all transformations of the form (\ref{uaction})
is isomorphic to the special Euclidean group $\mathbb{S}\mathbb{E}(2)$ of all
planar translations and rotations.

When studying the effects of inhomogeneities on the propagation of spiral
waves, one must consider a larger class of models than (\ref{basicrdpde}),
since inhomogeneous media do not possess Euclidean invariance. For example, one
might consider systems of the form
\begin{equation}
\frac{\partial u}{\partial
t}=D\cdot\nabla^2\,u+f(u)+\lambda\,g(u,\|x\|^2,\lambda) \label{basicrdpdep1}
\end{equation}
which are perturbations of (\ref{basicrdpde}).  Such systems could model a
spatially extended reaction-diffusion medium in which there is one site of
inhomogeneity (with circular symmetry) centered at the origin of
$\mathbb{R}^2$.  For instance, the Oregonator model which is used to study
spiral anchoring in \cite{MPMPV} is of the form (\ref{basicrdpdep1}). When
$\lambda\neq 0$, (\ref{basicrdpdep1}) has rotational symmetry about the origin,
but does not possess any translation symmetry.  This phenomenon is called {\it
forced translational symmetry-breaking}; it is studied in detail in \cite{LW}.
\par In this paper, we use a similar equivariant dynamical systems approach to study the
problem of spiral wave dynamics (specifically, with regards to
anchoring/repelling) in media in which there are several sites of
inhomogeneities (as opposed to just one site), of which cardiac tissue is an
important example.

We will make several simplifying assumptions which are meant to make the
analysis more tractable. First, we assume that the inhomogeneities consist of a
finite number of ``sources'' which are localized near distinct sites
$\zeta_1,\ldots, \zeta_n$ in the plane. Second, we will assume that these $n$
sources of inhomogeneity are independent in the following sense: we introduce
$n$ independent real parameters $\lambda_1,\ldots,\lambda_n$ which give some
measure of the relative ``amplitudes'' of the sources. In particular, when all
the $\lambda_i$ are zero except, say $\lambda_{i*}\neq 0$, then there is only
one source of inhomogeneity localized near the point $\zeta_{i*}$. In that
case, we will make a third simplifying assumption: the single inhomogeneity is
circularly symmetric around the point~$\zeta_{i*}$.  The following is an
example of a class of reaction-diffusion partial differential equations which
are perturbations of (\ref{basicrdpde}) and which might model such a situation:
\begin{equation}\label{basicrdpdep2}
u_t=\tilde{D}\Delta
u+f(u)+\sum_{j=1}^n\lambda_j\left[\hat{D}_j(\|x-\zeta_j\|^2,\lambda)\Delta
u+f_j(u,\|x-\zeta_j\|^2,\lambda)\right],\end{equation} where the functions
$\hat{D}_j,f_j$ are bounded and smooth enough. The goal of this paper is to
provide a detailed analysis of a larger class of abstract dynamical systems
which share the symmetry properties of (\ref{basicrdpdep2}):
\begin{enumerate}[(S1)]
\item when $\lambda_1=\cdots=\lambda_n=0$, the systems are
invariant under the action (\ref{uaction}) of the group
$\mathbb{S}\mathbb{E}(2)$,
\item when all the $\lambda_i$ are zero except $\lambda_{i*}$, the
systems have rotational symmetry about the point $\zeta_{i*}$, but they do not
generically possess translation symmetries,
\item when two or more of the $\lambda_i$ are non-zero, the
systems do not generically possess any of the symmetries (\ref{uaction}) except
for the identity.
\end{enumerate}
Our results will apply to the subclass $\mathsf{LC}_0$ of systems whose members
also generate a smooth local semi-flow on a suitable function space
\cite{SSW1,SSW2,SSW3,SSW4,FSSW}, as well as some technical conditions which
will be specified as we proceed.

The paper is organized as follows. In the second section, we derive the center
bundle equations of the semi-flow of a system in $\mathsf{LC}_0$, near a
hyperbolic rotating wave. We state and prove our main results in the third
section: to wit, spiral anchoring is generic in a parameter wedge. Then, we
provide a visual criterion characterizing the anchoring wedges in the case
$n=2$. Finally, we perform numerical experiments demonstrating the validity of
our results.


\normalsize\section{Reduction to the Center Bundle Equations}

Let $X$ be a Banach space, ${\mathcal U}\subset\mathbb{R}^n$ a neighborhood of
the origin and $\Phi_{t,\lambda}$ be a smoothly parameterized family
(parameterized by $\lambda\in {\mathcal U}$) of smooth local semi-flows on~$X$.

Let $\mathbb{S}\mathbb{E}(2)=\mathbb{C}\dot{+}\mathbb{S}\mathbb{O}(2)$ denote
the group of all planar translations and rotations, and let
\begin{equation}
a:\mathbb{S}\mathbb{E}(2)\longrightarrow \GL(X) \label{a_action}
\end{equation}
be a faithful and isometric representation of $\mathbb{S}\mathbb{E}(2)$ in the
space of bounded, invertible linear operators on $X$.  For example, if $X$ is a
space of functions with planar domain, a typical $\mathbb{S}\mathbb{E}(2)$
action (such as (\ref{uaction}) in the preceding section) is given by
$$
(a(\gamma)u)(x)=u(\gamma^{-1}(x)),\,\,\,\,\gamma\in\mathbb{S}\mathbb{E}(2).
$$
We will parameterize $\mathbb{S}\mathbb{E}(2)$ as follows:
$\mathbb{S}\mathbb{E}(2)\cong\,\mathbb{C}\times\mathbb{S}^1$, with
multiplication given by $(p_1,\varphi_1)\cdot (p_2,\varphi_2)=(
e^{i\varphi_1}p_2+p_1,\varphi_1+\varphi_2)$,
$\forall\,(p_1,\varphi_1),\,(p_2,\varphi_2)\in\,\mathbb{C}\times\mathbb{S}^1$.
For fixed $\xi\in\mathbb{C}$, we define the following subgroup of
$\mathbb{S}\mathbb{E}(2)$:
$$
\mathbb{S}\mathbb{O}(2)_{\xi}=\{\,(\xi,0)\cdot (0,\theta)\cdot
(-\xi,0)\,\,|\,\,\theta \in\mathbb{S}^1\,\}
$$
which is isomorphic to $\mathbb{S}\mathbb{O}(2)$, and represents rotations
about the point $\xi$.  We will assume the following symmetry conditions on the
family $\Phi_{t,\lambda}$ of semi-flows.
\begin{hyp} \label{hyp1}
There exists $n$ distinct points $\xi_1,\ldots,\xi_n$ in $\mathbb{C}$ such that
if $e_j$ denotes the $j^{\mbox{\footnotesize th}}$ vector of the canonical
basis in $\mathbb{R}^n$, then $\forall\,u\in\,X, \alpha\neq 0, t>0,$
\begin{align*}
\Phi_{t,\alpha e_j}(a(\gamma)u)&=a(\gamma)\Phi_{t,\alpha
e_j}(u) \iff \gamma\in\,\mathbb{S}\mathbb{O}(2)_{\xi_j}, \quad\mbox{and}\\
\Phi_{t,0}(a(\gamma)u)&=a(\gamma)\Phi_{t,0}(u), \quad \forall\,
\gamma\in\,\mathbb{S}\mathbb{E}(2).
\end{align*}
\end{hyp}
Hypothesis \ref{hyp1} basically states that (a) when $\lambda=0$, the semi-flow
$\Phi_{t,0}$ is $\mathbb{S}\mathbb{E}(2)$-equivariant; (b) when $\lambda\neq 0$
is near the origin and along the $j^{\mbox{\footnotesize th}}$ coordinate axis
of $\mathbb{R}^n$,  the semi-flow is only
$\mathbb{S}\mathbb{O}(2)_{\xi_j}$-equivariant (i.e. it only commutes with
rotations about the point $\xi_j$), and (c) when $\lambda$ is not as in (a) or
(b), the semi-flow has (generically) trivial equivariance.

We are interested in the effects of the forced symmetry-breaking on normally
hyperbolic rotating waves.  Therefore, we will assume the following hypothesis.
\begin{hyp}
There exists $u^*\in X$ and $\Omega^*$ in the Lie algebra of
$\mathbb{S}\mathbb{E}(2)$ such that $ e^{\Omega^*t}$ is a rotation and
$\Phi_{t,0}(u^*)=a( e^{\Omega^*t})u^*$ for all $t$. We also assume that the set
$\{\,\lambda\in\mathbb{C}\,\,|\,\,|\lambda|\geq 1\,\}$ is a spectral set for
the linearization $a( e^{-\Omega^*})D\Phi_{1,0}(u^*)$ with projection $P_*$
such that the generalized eigenspace $\mbox{\rm range}(P_*)$ is three
dimensional. \label{hyp2}
\end{hyp}

For sake of simplicity, we will only be interested in one-armed spiral waves;
therefore, we assume that $u^*$ in hypothesis \ref{hyp2} has trivial isotropy
subgroup. While hypotheses 1 and 2 hold for a large variety of spirals (such as
decaying spirals), there is also a large family of spirals for which they don't
(including Archimedean spirals) \cite{S}.\footnote{However, even in the case of
Archimedean spirals (for which hypothesis 1 fails), finite-dimensional
center-bundle equations which share the symmetries of the underlying abstract
dynamical systems have been shown to possess a definite predictive value in
terms of possible dynamics and bifurcations of these spiral waves
\cite{B1,B2,LOPS,LW,LeB}.} \newl Let $\mathsf{LC}_0$ be the collection of all
abstract dynamical systems that do satisfy them, as well as all other
hypotheses required in order for the center manifold theorems of
\cite{SSW1,SSW2,SSW3,SSW4} to hold, and let
 $\Phi_{t,\lambda}$ be produced by some member of $\mathsf{LC}_0$. It follows that for $\lambda$ near
the origin in $\mathbb{R}^n$, the essential dynamics of the semi-flow
$\Phi_{t,\lambda}$ near the rotating wave reduces to the following ordinary
differential equations on the bundle $\mathbb{C}\times\mathbb{S}^1$ (see
\cite{Byeah} for more details):
\begin{align}
\begin{split} \label{basiceqs1}
\dot{p}&= e^{i\varphi}V+G^p(p,\overline{p},\varphi,\lambda)\\
\dot{\varphi}&=\omega+G^{\varphi}(p,\overline{p},\varphi,\lambda)
\end{split}
\end{align}
where $V$ is a complex constant, $\omega\neq 0$ is a real constant, $G^p$ and
$G^{\varphi}$ are smooth, uniformly bounded in $p$, and such that
$G^p(p,\overline{p},\varphi,0)=0$ and
$G^{\varphi}(p,\overline{p},\varphi,0)=0$.  If $\lambda$ is near the origin, we
can re-scale time along orbits of (\ref{basiceqs1}) to get
\begin{align}
\begin{split}\label{basiceqs2}
\dot{p}&= e^{i\varphi}v+{\mathcal G}(p,\overline{p},\varphi,\lambda)\\
\dot{\varphi}&=1
\end{split}
\end{align}
where ${\mathcal G}$ is smooth, uniformly bounded in $p$, and such that
${\mathcal G}(p,\overline{p},\varphi,0)\equiv 0$. Of course, ${\mathcal G}$ is
not completely arbitrary because of the symmetry conditions in hypothesis
\ref{hyp1}. A simple computation and Taylor's theorem lead to the following.
\begin{prop}
The symmetry conditions in hypothesis $\ref{hyp1}$ imply that the equations
$(\ref{basiceqs2})$ have the general form
\begin{align}
\label{system1} \dot{p}&= e^{i\varphi(t)}\left[v+
\sum_{j=1}^n\,\lambda_jH_j((p-\xi_j) e^{-i\varphi(t)},(\overline{p-\xi_j})
e^{i\varphi(t)}, \lambda)\right]
\end{align}
where, without loss of generality, $\varphi(t)=t$, $v\in \GC$,
$\lambda=(\lambda_1,\ldots,\lambda_n)$, and the functions $H_j$ are smooth and
uniformly bounded in $p$.
\end{prop}
A $2\pi-$periodic solution $p_{\lambda}$ of $(\ref{system1})$ is called a
\textit{perturbed rotating wave} of $(\ref{system1})$. Define the average value
\begin{align}[p_{\lambda}]_{\A}&=\frac{1}{2\pi}\int_{0}^{2\pi}\!\!\!\!p_{\lambda}(t)\, dt.\label{ca}\end{align}
If the Floquet multipliers of $p_{\lambda}$ all lie within (resp. outside) the
unit circle, we shall say that $[p_{\lambda}]_{\A}$ is the \textit{anchoring}
(resp. \textit{repelling}, or \textit{unstable anchoring}) \textit{center}
of~$p_{\lambda}$.\par In the following section, we will perform an analysis of
anchoring of perturbed rotating waves of (\ref{system1}) for parameter values
near $\lambda=0$.


\normalsize\section{Analysis of the Center Bundle Equations}

Equations (\ref{system1}) represent the dynamics near a normally hyperbolic
rotating wave for a parameterized family $\Phi_{t,\lambda}$ of semi-flows
satisfying the forced-symmetry breaking conditions in hypothesis \ref{hyp1}. We
start with a brief review of the case $n=1$ which was studied in detail in
\cite{LW}, and then present new results on the general $n$ case.

\normalsize\subsection{The Case $n=1$}

In this case, we may assume without loss of generality that $\xi_1=0$, so that
(\ref{system1}) has the form
\begin{align}
\label{basiceqs4} \dot{p}&={\displaystyle  e^{i t}\left[v+ \lambda H(p e^{-i
t},\overline{p} e^{i t}, \lambda)\right]}
\end{align}
where $\lambda\in\mathbb{R}$ is small.  By writing $w=p e^{-i t}\!\!+i v$, this
system becomes
\begin{align}
\label{basiceqs5} \dot{w}&={\displaystyle -i
w+\lambda{\widetilde{H}}(w,\overline{w},\lambda)}
\end{align}
where $\widetilde{H}(w,\overline{w},\lambda)= H(w-i
v,\overline{w}+i\overline{v},\lambda)$. The following theorem is proved in
\cite{LW}.
\begin{thm}
Let $a=\mbox{\rm Re}(D_1\widetilde{H}(0,0,0))$, where $\tilde{H}$ is as in
$(\ref{basiceqs5})$.  If $a\neq 0$, then for all $\lambda\neq 0$ small enough,
 $(\ref{basiceqs4})$ has a hyperbolic rotating wave
\begin{equation}
p(t)=\left(-i v+O(\lambda)\right) e^{i t},\quad\varphi(t)=t.
\label{anchored_sol}
\end{equation} The origin $[p]_{\A}=0$ is an anchoring center if $a\lambda<0$; it is a
repelling center if $a\lambda>0$.
\end{thm}
\begin{rem}
In the case where the semi-flow $\Phi_{t,\lambda}$ is generated by a system of
planar reaction-diffusion partial differential equations, the solution
(\ref{anchored_sol}) represents a wave which is rigidly and uniformly rotating
around the origin in the plane. In the case where $a\lambda<0$, the rotating
wave is locally asymptotically stable. When $a\lambda>0$, the rotating wave is
unstable (see \cite{MPMPV} for an experimental characterization of this
phenomenon in an Oregonator model).
\end{rem}

\normalsize\subsection{The Case $n>1$}

One might think that the combination of many perturbations would just combine
the effects of each perturbation, so that spirals would be observed anchoring
at each of the centers, but we shall see that this is not usually the case.\par
By re-labeling the indices in (\ref{system1}) if necessary, we can temporarily
shift our point of view so that $\xi_1$ plays the central role in the following
analysis. Then, under the co-rotating frame of reference $z=p-\xi_1+i e^{i
t}v$, (\ref{system1}) becomes
\begin{align}\label{zdotforced2rescaled}
\dot{z}=\dot{p}- e^{i t}v= e^{i t}\sum_{j=1}^n \lambda_j H_j \big((z-\zeta_j)
e^{-i t}\!\!-i v,\ov{(z-\zeta_j)} e^{i t}\!\!+i\ov{v},\lambda\big),
\end{align}
where $\zeta_j=\xi_j-\xi_1$ for $j=1,\ldots, n$. \newl When $\lambda_1\neq 0$
and $\lambda_2=\cdots=\lambda_n=0$, we find ourselves in the situation
described in the previous subsection. Now, set $\varepsilon=\lambda_1$,
$\mu_1=1$ and $\lambda_j=\mu_j\varepsilon$ for $j=2,\ldots, n$ and
$\mu=(\mu_2,\ldots,\mu_n)\in \GR^{n-1}$. Then (\ref{zdotforced2rescaled}) can
be viewed as a perturbation of the corresponding equation in the case $n=1$.
Note that $\zeta_1=0$ and
 $\lambda=(1,\mu)\varepsilon$. \newl Equation
(\ref{zdotforced2rescaled}) rewrites as
\begin{equation}\label{system2}
\dot{z}=\varepsilon  e^{i t}\sum_{j=1}^n \mu_j H_j \big((z-\zeta_j) e^{-i
t}\!\!-i v,\ov{(z-\zeta_j)} e^{i t}\!\!+i\ov{v},(1,\mu)\varepsilon\big).
\end{equation}
Let $\label{funcH} \hat{H}_j(w,\ov{w},\varepsilon,\mu)=H_j\big(w-i
v,\ov{w}+i\ov{v},(1,\mu)\varepsilon\big)$ for $j=1,\ldots n$. Then
(\ref{system2}) becomes
\begin{align}\label{system3}
\dot{z}=\varepsilon  e^{i t}K(z e^{-i t},\ov{z} e^{i t},t,\varepsilon,\mu)
\end{align}
where $\displaystyle{K(w,\ov{w},t,\varepsilon,\mu)=\sum_{j=1}^n \mu_j
\hat{H}_{j}(w-\zeta_j e^{-i t},\ov{w}-\ov{\zeta}_j e^{i t},\varepsilon,\mu)}$
is $2\pi-$periodic in $t$.
\newl Set $\alpha_1=D_1H_1(-i v,i\ov{v},0)$. The time$-2\pi$ map $P$ of
(\ref{system3}) is given by
\begin{equation}\label{Poincare2}
P(z,\overline{z},\varepsilon,\mu)=z+2\pi\varepsilon\Big[\alpha_1z+O\big(|z|^2\big)+O\big(\varepsilon,\mu_2,\ldots,\mu_n\big)\Big]
\end{equation} near $z=0$ and $(\varepsilon,\mu)=(0,0)$. \par Hyperbolic fixed points of (\ref{Poincare2})
correspond to hyperbolic $2\pi-$periodic solutions of (\ref{system3}), and so
to perturbed rotating waves of (\ref{system1}), that is, the path traced by the
solution wave need not be circular. As $z=0$ is not generally a fixed point of
(\ref{Poincare2}), these perturbed rotating waves may not be centered
at~$\xi_1$. Indeed, let
\begin{equation}\label{funcB}B(z,\overline{z},\varepsilon,\mu)=\alpha_1z+O\big(|z|^2\big)+O\big(\varepsilon,\mu_2,\ldots,\mu_n\big)\end{equation} be the function inside the square brackets in (\ref{Poincare2}). Note that
$B(0,0,0,0)=0$ and that, generically, $D_1B(0,0,0,0)=\alpha_1\neq 0$. By the
implicit function theorem, there is a unique smooth function
$z(\varepsilon,\mu)$ defined near $(\varepsilon,\mu)=(0,0)$ with $z(0,0)=0$ and
\begin{align}\label{funcBB}B\big(z(\varepsilon,\mu),\overline{z}(\varepsilon,\mu),\varepsilon,\mu\big)\equiv
0\end{align} near $z=0$. This leads to the following theorem.
 \begin{thm} \label{thm41}Let $\alpha_1$ be as in the preceding discussion,
   with $\mbox{\rm Re}(\alpha_1)\neq 0.$ If the parameters are small
   enough to satisfy the conditions outlined in the proof below, the time$-2\pi$ map $(\ref{Poincare2})$ has a unique family of hyperbolic fixed points, whose stability is exactly determined by the sign of $\varepsilon \real(\alpha_1)$.
\end{thm}
  \noindent\proof Let $B$ be as in
(\ref{funcB}) and $z(\varepsilon,\mu)$ be the unique continuous function
solving the equation $B=0$ for small parameter values, as asserted above. When
$\varepsilon=0$, any point of $\GR^2$ is a non-hyperbolic fixed point of $P$
and so, from now on, we will assume that $\varepsilon\neq 0$. If that is the
case, and if $\varepsilon$ and $\|\mu\|$ are small enough, the eigenvalues
$\omega_{1,2}(\varepsilon,\mu)$ of $DP(z(\varepsilon,\mu),\varepsilon,\mu)$
satisfy
\begin{align*}
\left|\omega_{1,2}(\varepsilon,\mu)
\right|^2&=1+4\pi\varepsilon\real(\alpha_1)+\varepsilon O(\varepsilon,\mu)\neq
1,
\end{align*} since $\real(\alpha_1)\neq 0$. In other words, the fixed point $z(\varepsilon,\mu)$ is hyperbolic. When
$\varepsilon\real(\alpha_1)<0$, the eigenvalues lie inside the unit circle and
the fixed point is asymptotically stable; otherwise, it is unstable. \qed
\newl We are now able to formulate and prove the following result.
\begin{thm} \label{thm42} Suppose the hypotheses of theorem~$\ref{thm41}$ are satisfied. Then there exists a wedge-shaped region near
    $\lambda=0$ of the form
$$
{\mathcal W}_{1}=\{(\lambda_1,\ldots,\lambda_n)\in
\GR^n\,:\,|\lambda_j|<W_{1,j}|\lambda_1|,\,\,\,W_{1,j}>0,\,\,\mbox{ \rm for
$j\neq 1$ and $\lambda_1$ near}\,\,0\,\}
$$
such that for all $0\neq \lambda\in {\mathcal W}_{1}$, $(\ref{system1})$ has a
unique perturbed rotating wave $\mathcal{S}^1_{\lambda}$, with center
$[\mathcal{S}^1_{\lambda}]_{\A}$ generically away from $\xi_1$. Furthermore,
$[\mathcal{S}^1_{\lambda}]_{\A}$ is a center of anchoring when
$\lambda_1\real(\alpha_1)< 0$.
\end{thm}
\noindent\proof For $j\neq 1$, let $W_{1,j}> 0$ be such that the conclusion of
theorem~\ref{thm41} holds for any $\mu_j$ with $|\mu_j|<W_{1,j}$. Let
$\mathcal{W}_1$ be as stated in the hypothesis. If $(\varepsilon,\mu)$ is such
that the time$-2\pi$ map (\ref{Poincare2}) has a hyperbolic fixed point
$z(\varepsilon,\mu)$ near $0$, then (\ref{system2}) has a hyperbolic
$2\pi-$periodic orbit $\tilde{z}_{\varepsilon,\mu}(t)$ centered at a point near
$z=\zeta_1=0$. \par\noindent For $j\neq 1$, let $\lambda_1=\varepsilon\neq 0$
be small enough  and set $\lambda_j=\mu_j\varepsilon$. Then $\lambda\in
\mathcal{W}_1$, as
$$|\lambda_j|=|\mu_j||\varepsilon|<W_{1,j}|\lambda_1|\quad \mbox{for
}j\neq 1,$$ and $\tilde{z}_{\varepsilon,\mu}(t)$ is a $2\pi-$periodic orbit for
the parameter $\lambda$, which we denote by $z_{\lambda}(t)$. Since $p=z-i e^{i
t}v+\xi_1$, (\ref{system1}) has a unique perturbed rotating wave
$\mathcal{S}_{\lambda}^1$, with
$$[\mathcal{S}_{\lambda}^1]_{\A}=\frac{1}{2\pi}\int_{0}^{2\pi}\!\!\!\!\left(z_{\lambda}(t)-i e^{i t}v+\xi_1\right)\,
dt=\xi_{1}+[z_{\lambda}]_{\A}.$$ If $0\neq\lambda\in \mathcal{W}$ is such that
$\mu_j=\lambda_j/\varepsilon\neq 0$ is fixed for $j=2,\ldots,n$, then
$[z_{\lambda}]_{\A}=O(1)$ as $\lambda_1\to 0$ and so
$[\mathcal{S}_{\lambda}^1]_{\A}\neq \xi_1$, generically. The conclusion about
the stability of ${\mathcal S}^1_{\lambda}$ follows directly from
theorem~\ref{thm41}. \qed \begin{rem} When $\lambda$ approaches the
$\lambda_1-$axis away from the origin, $[\mathcal{S}^1_{\lambda}]_{\A}\to
\xi_1$. On the other hand, when the parameter values stray outside of
${\mathcal W}_{1}$, all that can generically be said with certainty is that
solutions of (\ref{system1}) locally drift away from $\xi_1$, which cannot then
be a center of anchoring. After drifting, the spiral may very well get anchored
at some point far from $\xi_1$, depending on the global nature of the
perturbation functions $H_j$ in (\ref{system1}).\end{rem} The preceding results
have been achieved by considering (\ref{system1}) under a co-rotating frame of
reference around $\xi_1$. Of course, since the choice for $\xi_1$ was
arbitrary, corresponding results must also be achieved, in exactly the same
manner, when the viewpoint shifts to another~$\xi_k$. For $j=1,\ldots, n,$ let
$\alpha_j=D_1H_j(-i v,i\ov{v},0)$ be the \textit{anchoring coefficients} of
(\ref{system1}).

\begin{thm}\label{thm43} Let $k\in \{1,\ldots,n\}$. If $\real(\alpha_k)\neq 0$,
 then there exists a wedge-shaped region near
    $\lambda=0$ of the form
$${\mathcal W}_{k}=\{(\lambda_1,\ldots,\lambda_n)\in
\GR^n\,:\,|\lambda_j|<W_{k,j}|\lambda_k|,\,\,\,W_{k,j}>0,\,\,\mbox{\rm for
$j\neq k$ and $\lambda_k$ near}\,\,0\,\}
$$
such that for all $0\neq \lambda\in {\mathcal W}_{k}$, $(\ref{system1})$ has a
unique perturbed rotating wave $\mathcal{S}^k_{\lambda}$, with center
$[\mathcal{S}^k_{\lambda}]_{\A}$ generically away from $\xi_k$. Furthermore,
$[\mathcal{S}^k_{\lambda}]_{\A}$ is a center of anchoring when
$\lambda_k\real(\alpha_k)< 0$.
\end{thm}
Clearly, the remark that appears after the proof of theorem~\ref{thm42} still
holds.

\normalsize


\section{Characterization of Spiral Anchoring $(n=2)$}\label{CSA}
\normalsize In the previous section, we described the (local) behaviour of
spiral anchoring in small wedges around the parameter coordinate axes. In this
section, we present a fuller characterization of spiral wave anchoring for the
case $n=2$.\footnote{Most of the analysis can be extended and adapted to the
general case $n\geq 2$, but at the cost of substantial algebraic
complications.}\newl Let $0\neq \xi\in\GR^2$,
$\Lambda_0=(\lambda_1,0),\Lambda_{\xi}=(0,\lambda_2)\in \GR^2$ and let
$P:\GR^2\times\GR^2\to \GR^2$ be a real analytic map with $P(x,0)=x$,
$DP(x,0)=I_2$ for all $x\in \GR^2$, satisfying the following conditions: for
$\eta\in\{0,\xi\}$,
\begin{enumerate}[(P1)]
\item $\exists\, \omega_*>0$ such that $P(\eta,\Lambda_{\eta})\equiv 0$, for all $||\Lambda_{\eta}||<
\omega_*$;
\item the eigenvalues of $DP(\eta,\Lambda_{\eta})$ lie both outside or both inside the unit circle
for all $0\neq ||\Lambda_{\eta}||< \omega_*$;
\item there is a wedge region $\mathsf{w}_{\eta}$
surrounding the coordinate axis generated by $\Lambda_{\eta}$ in parameter
space \begin{figure}[t]
\begin{center}
\includegraphics[width=175pt]{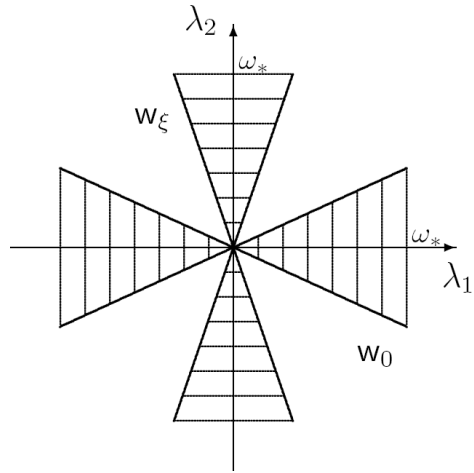}
\caption{ Wedges in parameter space corresponding to property
(P3).}\end{center}\hrule
\end{figure}\normalsize(see figure~\thefigure) in which $P$ has a (locally) unique
manifold $x_{\eta}(\lambda)$ such that, for all $\lambda\in \mathsf{w}_{\eta}$,
\begin{enumerate}
\item $P(x_{\eta}(\lambda),\lambda)\equiv
x_{\eta}(\lambda)$;
\item $x_{\eta}(\lambda)\to \eta$ as $\lambda$ approaches the coordinate axis away from the origin;
\item $x_{\eta}(\lambda)$ shares its stability with $\eta$ in (P2).
\end{enumerate}
\end{enumerate}\normalsize
When the hypotheses of theorem \ref{thm43} hold, the associated time$-2\pi$ map
(\ref{Poincare2}) (viewed in real coordinates) satisfies (P1)$-$(P3). Numerous
questions cannot be answered by local analysis alone. For instance:
\begin{enumerate}
\item Can the wedges overlap? What does that imply for anchoring in
(\ref{system1})?
\item Can a wedge contain its ``opposite'' coordinate axis?
\item If the wedges do not overlap, what is the nature of their complement?
\item If there is a complement with non-trivial measure, what
kinds of dynamics can be expected as the parameter vector $\lambda$ traces a
circle around the origin in parameter space?
\end{enumerate}
We will provide answers to these questions by first studying a specific map,
then extending our results to the general mapping. \normalsize\subsection{A
Specific Mapping}\label{ASM} \normalsize\sectionmark{A Specific
Mapping}\normalsize Consider the mapping
$P:\GR^2\times\GR^2\to\GR^2$\index{P@$P$} given by
\begin{equation}\label{themapping}
P(x,\lambda)=x+2\pi\big[\lambda_1F_0(x)+\lambda_2G_{\xi}(x)\big],
\end{equation} where $0\neq \xi\in \GR^2$, and $F_0$, $G_{\xi}$
are real analytic functions of $x,\lambda\in \GR^2$. \par Such a map is
obtained by truncating the $\lambda-$terms of order $\geq 2$ from the
time$-2\pi$ map (\ref{Poincare2}), for instance. According to
theorem~\ref{thm43}, the jacobians $DF_0(0)$ and $DG_{\xi}(\xi)$ have a
particular structure. \normalsize
\begin{prop}\label{themappingprop}
If $F_0(0)=0$, $G_{\xi}(\xi)=0$, and if $$DF_0(0)=\begin{pmatrix}a & -b
\\ b & a\end{pmatrix}\quad\mbox{and}\quad
DG_{\xi}(\xi)=\begin{pmatrix}c & -d \\ d & c\end{pmatrix}$$ where $a,c\neq 0$,
then there exists $\omega_*>0$ such that the map defined by
$(\ref{themapping})$ satisfies the conditions
\textsc{(P1)$-$(P3)}.\end{prop}\normalsize \normalsize\subsubsection{The search
for fixed points.}\label{TSFFP} Define $A:\GR^2\to
 \GM_2(\GR)$ by \begin{equation} A(x)=\begin{bmatrix} F_0(x) &
 G_{\xi}(x)\end{bmatrix}.\label{theA}\end{equation} Then, $\hat{x}$
 is a fixed point of (\ref{themapping}) for $\hat{\lambda}\in \GR^2$
 if and only if $A(\hat{x})\cdot{\hat{\lambda}}=0$, that is if and only if
 $\hat{\lambda}\in L_{\hat{x}}=\ker A(\hat{x}).$\index{L-x@$L_{\hat{x}}$}
Let $(\hat{x},\hat{\lambda})$ be such a pair. According to the implicit
function theorem, as long as \begin{align}\label{theimplicit}\det
\left(D_xP(\hat{x},\hat{\lambda})-I\right)=4\pi^2\det
\left(\hat{\lambda}_1DF_0(\hat{x})+\hat{\lambda}_2DG_{\xi}(\hat{x})\right)\neq
0,\end{align} there is a neighbourhood $W$ of $\hat{\lambda}$ and a unique
analytic function $X:W\to \GR^2$ such that $X(\hat{\lambda})=\hat{x}$ and
$A(X(\lambda))\cdot \lambda \equiv 0 \mbox{ for all }\lambda\in W.$ By
construction, $X(\lambda)$ is a fixed point of (\ref{themapping}) for all
$\lambda\in W$.\par If $\dim L_{\hat{x}}=0$ as a manifold, then
$L_{\hat{x}}=\{0\}$. Consequently, the preceding implicit function theorem
construction fails, which contradicts property (P3). We need thus only
investigate fixed points $\hat{x}$ for which $\dim L_{\hat{x}}\neq 0$. As the
quantities under consideration are analytic, it can further be assumed that
$\rank A(\hat{x})=1$ and $\dim L_{\hat{x}}=1$.
\newl We now show how to optimally extend the wedge regions $\mathsf{w}_{\eta}$ using
property (P3). Let $(x^*,\lambda_*),(x_*,\lambda_*)\in \GR^2\times
(\GR^2-\{0\})$ be such that $x^*$, $x_*$ are fixed points of
(\ref{themapping}), $\lambda^*\in L_{x^*},\lambda_*\in L_{x_*}$, and
(\ref{theimplicit}) is satisfied for both pairs.  According to the implicit
function theorem, there are open neighbourhoods $W^*,W_*$ of
$\lambda^*,\lambda_*\in\GR^2$ respectively, and a pair of unique real analytic
functions $X^*:W^*\to \GR^2$, $X_*:W_*\to \GR^2$ for which
$X^*(\lambda^*)=x^*$, $X_*(\lambda_*)=x_*$ and
\begin{align*}A(X^*(\Lambda))\cdot \Lambda &\equiv 0,\ \mbox{for all }\Lambda\in
W^*,\quad  A(X_*(\Lambda))\cdot \Lambda \equiv 0,\ \mbox{for all }\Lambda\in
W_*. \end{align*}
\begin{lemma}\label{thelemma} If $\Lambda_*^*\in W^*_*=W^*\cap W_*$ is such that
$X^*(\Lambda_*^*)=X_*(\Lambda_*^*)$, then $X^*=X_*$ on~$W^*_*$
\end{lemma}
\noindent\proof The assertion follows from the uniqueness of the real analytic
functions $X^*, X_*$ in the implicit function theorem.\qed\normalsize  \newl
Denote the punctured open disc of radius $\omega_*$ centered at the origin by
$B(0,\omega_*)$. Let $\eta\in \{0,\xi\}$, $\omega_*>0$ be as in
proposition~\ref{themappingprop} and $0\neq \Lambda_{\eta}\in
\mathsf{w}_{\eta}$ be a point on the appropriate coordinate axis, as in
properties (P1) and (P2). According to these same properties, $\eta$ is a fixed
point of (\ref{themapping}) for $\Lambda_{\eta}$ and $\det
\left(D_xP(\eta,\Lambda_{\eta})-I\right) \neq 0.$\par Lemma~\ref{thelemma} then
implies the existence of a maximal open region $\mathsf{W}_{\eta}$, defined as
a union of open sets $W\subseteq B(0,\omega_*)$ (in much the same way as the
maximal interval is built in the Fundamental Theorem of ODE \cite{HS}),
containing $\mathsf{w}_{\eta}\cap B(0,\omega_*)$ and for which there is a
unique real analytic function $\mathsf{X}_{\eta}:\mathsf{W}_{\eta}\to \GR^2$
satisfying $x_{\eta}=\mathsf{X}_{\eta}|_{\mathsf{W}_{\eta}}$, where $x_{\eta}$
is as in property (P3). \begin{figure}[t]
\begin{center}
\includegraphics[width=150pt]{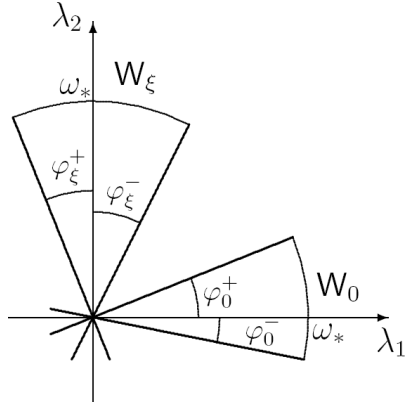}
\caption{ Wedge angles, with optimal wedge-like regions in parameter
space.}\end{center}\hrule
\end{figure}\noindent Since $\hat{x}$ is a fixed point of (\ref{themapping}) for $0\neq \hat{\lambda}$
whenever $\hat{\lambda}\in L_{\hat{x}}$, $\mathsf{W}_{\eta}$ is described (in
polar coordinates) by either one of
\begin{align*}\mathsf{W}_{\eta}&=\{(r,\theta):0<r<\omega_*\mbox{ and }
s_{\eta}-\varphi_{\eta}^-<\theta<s_{\eta}+\varphi_{\eta}^+\} \\
\mathsf{W}_{\eta}&=\{(r,\theta):0<r<\omega_*\mbox{ and }\theta\in [0,2\pi]\}
 \end{align*} where $\varphi_{\eta}^-,\varphi_{\eta}^+\in (0,\pi/2]$ and
 \begin{align}s_{\eta}=\begin{cases}0 &\text{if $\eta=0$,} \\ \pi/2 &\text{if $\eta=\xi$.} \end{cases}\label{theseta},\end{align}
In the latter case, we will say that $\mathsf{W}_{\eta}$ is
\textit{catastrophe-free}. In the former case, the quantities
$\varphi_{\eta}^-,\varphi_{\eta}^+\in (0,\pi/2]$ are called the
\textit{fore-angle} and \textit{post-angle} of $W_{\eta}$, respectively (see
figure~\thefigure).\newl The implicit function theorem fails to extend
$A(\mathsf{X}_{\eta}(\lambda))\cdot \lambda \equiv 0$ (that is, it fails to
extend $\mathsf{W}_{\eta}$) at $(x^*,\lambda^*)$ if either
\begin{enumerate}[(C1)]
\item $\det \left(\lambda^*_1DF_0(x^*)+\lambda^*_2DG_{\xi}(x^*)\right)= 0$ and $\mathsf{X}_{\eta}(\lambda)\to x^*$ as $\lambda\to\lambda^*$, or
\item $\|\mathsf{X}_{\eta}(\lambda)\|\to\infty$ as $\lambda\to\lambda^*$.
\end{enumerate}
Such events will be referred to as \textit{fold} and
\textit{$\infty-$catastrophes}, respectively, or \textit{catastrophes},
collectively.
\par Let $0<\rho<\omega_*$ and set
\begin{align}\label{thecircle}\gamma_{\rho}(s)=\rho\left(
\cos(s),\sin(s) \right)^{\!\top}.\index{g-rho@$\gamma_{\rho}$}\end{align}
Assume $\mathsf{W}_{\eta}$ is not catastrophe-free. Starting at
$(\rho,s_{\eta})\in \mathsf{W}_{\eta}$, denote the angles in $(0,\pi/2]$
measuring the first clockwise and the first counter-clockwise occurrence of a
catastrophe along $\gamma_{\rho}$ by $\theta_{\eta}^-$ and $\theta_{\eta}^+$
respectively. Then, $\varphi_{\eta}^{\pm}=s_{\eta}\pm \theta_{\eta}^{\pm}$.
\normalsize\subsubsection{Fold bifurcation points.} Modulo a simple regularity
condition (see below), (C1) is equivalent to the existence of a fold
bifurcation curve in parameter space for (\ref{themapping}). Indeed, in that
case, $(x^*,\lambda^*)$ is a solution of
\begin{align}P(x,\lambda)-x&=0,\qquad  \det \left(D_xP(x,\lambda)-I\right)=0. \label{thefold}\end{align}
If the (full) Jacobian of the left-hand side of (\ref{thefold}) has rank 3 at
that point, (\ref{thefold}) has a fold bifurcation curve through $\lambda^*$
\cite{KU}. Such solutions are in one-to-one correspondence with regular
solutions of \begin{align}A(x)\cdot \lambda&=0,\qquad \det
\left(D_x\left[A(x)\cdot \lambda\right]\right)=0. \label{thefold2}\end{align}
Set
$$I_{10}=\begin{pmatrix}1 & 0 \\ 0 & 0\end{pmatrix},\quad
I_{01}=\begin{pmatrix}0 & 0 \\ 0 & 1\end{pmatrix},\quad
\hat{I}=\begin{pmatrix}0 & 1 \\ 1 & 0\end{pmatrix},\quad e_1=\begin{pmatrix}1
\\ 0\end{pmatrix}\quad\mbox{and}\quad e_2=\begin{pmatrix}0 \\ 1\end{pmatrix},$$
and define $H_1,H_2:\GR^2\to \GR^2$ by
\begin{align*}
H_1(x)&=\left[I_{10}A(x)+I_{01}A(x)\hat{I}\right]e_1,\qquad
H_2(x)=\left[I_{10}A(x)+I_{01}A(x)\hat{I}\right]e_2.
\end{align*} A quick computation \label{pdefold} shows that (\ref{thefold2}) can be written as
\begin{align}\label{thefold3}A(x)\cdot \lambda&=0,\qquad \lambda^{\!\top}Q(x)\lambda=0, \end{align} where \begin{align*}Q(x)=\begin{pmatrix}B(x) & \frac{1}{2}C(x) \\ \frac{1}{2}C(x) & E(x)\end{pmatrix} \end{align*} and
\begin{align}\label{theBCE}
\begin{split}
B(x)&=\det DF_0(x) \\
C(x)&= \det DH_1(x)+\det DH_2(x) \\
E(x)&=\det DG_{\xi}(x).
\end{split}\end{align}
Let $x^*$ be a fixed point of (\ref{themapping}) and denote
$K_{x}=\{\lambda:\lambda^{\!\top} Q(x)\lambda=0 \}$. Generically, $K_{x^*}$
consists of a single line or a pair of intersecting lines through the origin in
parameter space. Writing $\mathsf{L}_{x}= L_{x}\cap B(0,\omega_*)$ and
$\mathsf{K}_{x}= K_{x}\cap B(0,\omega_*)$, we can summarize the situation with
the following proposition.
\begin{prop}\label{theotherprop} If $(x^*,\lambda^*)$ is a regular solution of $(\ref{thefold2})$
with $\{0\}\neq \mathsf{L}_{x^*}\subseteq \mathsf{K}_{x^*}$, then
$(x^*,\lambda)$ is a fold bifurcation point of $(\ref{themapping})$ for all
$\lambda\in \mathsf{L}_{x^*}$.\end{prop} \normalsize\subsection{The Visual
Criterion}\label{TVC} Set $\mathfrak{Z}=\{(x,\lambda): P(x,\lambda)=x\mbox{ and
}\lambda\in\mathsf{L}_{x}\neq \{0\}\}$ and
$$\kappa(\mathfrak{Z})=\{x:\exists\,  \lambda\neq 0
\mbox{ such that }(x,\lambda)\in \mathfrak{Z}\}.$$ By construction,
$\kappa(\mathfrak{Z})$ is the zero-set of $\det A(x)$ in $\GR^2$ and $0,\xi\in
\kappa(\mathfrak{Z})$. Generically, $\kappa(\mathfrak{Z})$ is a collection
$\mathcal{C}$ of isolated planar curves, whose constituents come in two
varieties: bounded or unbounded.\footnote{Indeed, were any such curves to
intersect at $x_*$, $P$ would undergo a transcritical bifurcation along
$\kappa^{-1}(x_*)$. Such bifurcations are not generically permitted by (C1) and
(C2).} Denote this partition by $\mathcal{C}=\mathcal{C}_{\BB}\sqcup
\mathcal{C}_{\infty}$ and let $C_0, C_{\xi}$ be the curves in $\mathcal{C}$ for
which $0\in C_0$ and $\xi \in C_{\xi}$.\newl Let $\gamma_{\rho}:[0,2\pi]\to
\GR^2$ be the circle of radius $\rho$ around the origin, parameterized as in
(\ref{thecircle}). For each $(x,\lambda)\in \mathfrak{Z}$, define
$\mathsf{L}^{\rho}_{x}$ and $\mathsf{K}^{\rho}_{x}$ as the intersection of that
circle with $\mathsf{L}_{x}$ and $\mathsf{K}_{x}$, respectively, and let
$P_{\rho}:\GR^2\times [0,2\pi]\to \GR^2$ be given by
\begin{equation}\label{themapping22}
P_{\rho}(x,s)=x+2\pi\rho\big[\cos(s)F_0(x)+\sin(s)G_{\xi}(x)\big].
\end{equation}
Then $\mathsf{L}_{x}^{\rho}$ consists of two antipodal points $\{\pm
\alpha_{x,\rho}\}$, and the fixed points $(x,s)$ of $P_{\rho}$ are in
one-to-one correspondence with the `lines' of fixed points $(x,\mathsf{L}_{x})$
of $P$ for which $\mathsf{L}_x\neq \{0\}$ (see
proposition~\ref{theotherprop}).\par Set $\mathfrak{Z}_{\rho}=\{(x,s)\in
\GR^2\times[0,2\pi]: P_{\rho}(x,s)=x\mbox{ and } \gamma_{\rho}(s)\in
\mathsf{L}_x^{\rho}\}$ and
$$\kappa_{\rho}(\mathfrak{Z}_{\rho})=\{x:(x,s)\in
\mathfrak{Z}_{\rho}\}.$$ By construction,
$\kappa_{\rho}(\mathfrak{Z}_{\rho})=\kappa(\mathfrak{Z}).$ Thus for each $C\in
\mathcal{C}$, $\kappa_{\rho}^{-1}(C)$ is a branch of fixed points in the
bifurcation diagram of  $P_{\rho}$. According to section~\ref{TSFFP}, the
converse also holds: each branch of fixed points in the bifurcation diagram of
$P_{\rho}$ projects down \textsl{via} $\kappa_{\rho}$ to a curve in
$\mathcal{C}$. \newl The existence and location of fold catastrophes cannot be
read directly from $\mathcal{C}$, but the next proposition remedies that
situation.\par Let $(x^*,\alpha)\in \GR^2\times (\GR^2-\{0\})$ be such that
$\|\alpha\|=\rho$, $\alpha\in \mathsf{L}^{\rho}_{x^*}\subseteq
\mathsf{K}^{\rho}_{x^*}$. Recall that $A(x^*)\neq 0$. Then, $\left(
A_{j,1}(x^*)\ A_{j,2}(x^*)\right)\neq 0$ for some $j\in \{1,2\}$ and
\begin{align*} L_{x^*}=\{\lambda:A_{j,1}(x^*)\lambda_1+A_{j,2}(x^*)\lambda_2=0\}.\end{align*}
The function $\Gamma_j:\GR^2\to \GR^2$ defined by
\begin{equation}\Gamma_j(x)=\big[A^2_{j,2}(x)B(x)-A_{j,1}(x)A_{j,2}(x)C(x)+A^2_{j,1}(x)E(x)\big],
\end{equation} where $B,C$ and $E$ are as in (\ref{theBCE}), is
called the \textit{$j-$fold bifurcation function} of (\ref{themapping}). Let
$\mathcal{R}_j$ be the zero-set of $\Gamma_j(x)$ in $\GR^2$. We shall say that
$x^*$ is a \textit{transverse intersection} of $\kappa(\mathfrak{Z})$ and
$\mathcal{R}_j$ if $\det A(x^*)=\Gamma_j(x^*)=0$ and
$$\rank D\begin{pmatrix}\det A(x^*) \\ \Gamma_j(x^*)\end{pmatrix}=2.$$
\begin{prop} \label{therealprop} Let $j\in \{1,2\}$. If $x^*$ is a
transverse intersection of $\kappa(\mathfrak{Z})$ and $\mathcal{R}_j$ such that
$\left( A_{j,1}(x^*)\ A_{j,2}(x^*)\right)\neq 0$ and either
\begin{enumerate}[$(1)$]
\item $B(x^*)=0$ and $A_{j,1}(x^*)C(x^*)-A_{j,2}(x^*)E(x^*)=0$ or
\item $B(x^*)\neq 0$ and $C(x^*)^2-4B(x^*)E(x^*)\geq 0$,
\end{enumerate} then $P_{\rho}$ undergoes a fold catastrophe at $(x^*,s^*)$ for all $s^*$ such that
$\gamma_{\rho}(s^*)=\pm {\alpha}_{x^*,\rho}$. \end{prop} \noindent\proof By
re-labeling the terms if necessary, we may assume $A_{j,1}\neq 0$. There are
then two possibilities.
\begin{enumerate}\item  If $B=0$ and $A_{j,1}C-A_{j,2}E=0$, then
$$K_{x^*}=\{\lambda:\lambda_1\lambda_2C+\lambda_2^2E=0\} =\{\lambda:\lambda_2=0 \mbox{ or } \lambda_1C+\lambda_2E=0\}.$$
\begin{enumerate}
\item If $\lambda_2=0$, then
$L_{x^*}=\{(\lambda_1,0): \lambda_1A_{j,1}=0\}=\{0\}$ since $A_{j,1}\neq 0$.
But this contradicts the assumption  $\dim L_{x^*}=1$.
\item If $\lambda_1C+\lambda_2E=0$, then $$\rank \begin{pmatrix}A_{j,1} & A_{j,2} \\ C & E\end{pmatrix}=1.$$
\end{enumerate}
\item If $B\neq 0$ and $C^2-4BE\geq 0$, then
$$K_{x^*}=\left\{\lambda:\lambda_1=\frac{-C\pm\sqrt{C^2-4BE}}{2B}\lambda_2\right\}.$$
In this case,
\begin{align*}
4B\Gamma_j&=\left(-2A_{j,2}B+A_{j,1}C\right)^2-A_{j,1}^2(C^2-4BE) \\ &=
\left(-2A_{j,2}B+A_{j,1}\left(C+\sqrt{C^2-4BE}\right)\right)\\
&\qquad\qquad\qquad\qquad \cdot
\left(-2A_{j,2}B+A_{j,1}\left(C-\sqrt{C^2-4BE}\right)\right)=0
\end{align*} and so $$-\frac{A_{j,2}}{A_{j,1}}=\frac{-C+\sqrt{C^2-4BE}}{2B}\quad\mbox{or}\quad -\frac{A_{j,2}}{A_{j,1}}=\frac{-C-\sqrt{C^2-4BE}}{2B}.$$
\end{enumerate} In either cases, $L_{x^*}$ is contained in $K_{x^*}$; thus $\{0\}\neq
\mathsf{L}_{x^*}\subseteq \mathsf{K}_{x^*}$ and $\{0\}\neq
\mathsf{L}_{x^*}^{\rho}\subseteq \mathsf{K}_{x^*}^{\rho}$.  As $x^*$ is a
transverse intersection of $\kappa(\mathfrak{Z})$ and $\mathcal{R}_j$, it is
also a regular solution of (\ref{thefold3}); $(x^*,\mathsf{L}_{x^*})$ then
consists of fold bifurcation points of (\ref{themapping}), according to
proposition~\ref{theotherprop}. The desired conclusion follows from $\{\pm
\alpha_{x^*,\rho}\}=\mathsf{L}_{x^*}^{\rho}=\mathsf{L}_{x^*}\cap \gamma_{\rho}$
and from the correspondence between fixed points of $P_{\rho}$ and `lines' of
fixed points of $P$. \qed \newl By construction, the bifurcation diagram of
$P_{\rho}$ is $2\pi-$periodic~in~$s$. Consequently, elements of
$\mathcal{C}_{\BB}$ must be (bounded) loops and elements of
$\mathcal{C}_{\infty}$ must give rise to two $\infty-$catastrophes. Moreover,
the number of fold catastrophes on any given $C\in \mathcal{C}_{\BB}$ cannot be
odd as $C$ could not be a loop were that the case. Finally, note that
catastrophes cannot occur at $0$ or $\xi$ as this would contradict (P2) and
(P3). \normalsize\subsection{The Bifurcation Diagrams} Let $C_0$, $C_{\xi}$,
$s_0$ and $s_{\xi}$ be as defined previously. Set $\eta\in \{0,\xi\}$. By
definition, $C_{\eta}$ goes through $\eta$ at $s=s_{\eta}$. By (P3), the
wedges' angles $\varphi_{\eta}^{\pm}$ lie in $(0,\pi)$ or $(0,\pi]$ (when they
exist), according to whether they record fold or $\infty-$catastrophes,
respectively. Set $\nu_1=\varphi_0^++\varphi_{\xi}^-$ and
$\nu_2=\varphi_0^-+\varphi_{\xi}^+$. Then $\mathsf{W}_0$ and $\mathsf{W}_{\xi}$
overlap
\begin{enumerate}
\item in all four quadrants if and only if $\nu_1,\nu_2>\pi/2$;
\item in the first and third quadrants if and only if $\nu_1>\pi/2$ and $\nu_2\leq\pi/2$, and
in the second and fourth quadrants if and only if $\nu_1\leq\pi/2$ and
$\nu_2>\pi/2$.
\end{enumerate} If $\nu_j=\pi/2$, the wedges do not overlap but their complement has zero measure in a
neighbourhood of the origin. When the wedge angles $\varphi_{\eta}^{\pm}$ do
not exist, $\mathsf{W}_{\eta}$ is a deleted neighbourhood of the origin in
parameter space. \normalsize\subsubsection{The case $C_0\neq C_{\xi}$} In this
instance, it is sufficient to understand the bifurcation diagrams along a
single curve: the full picture can then be obtained by combining the diagrams
corresponding to $C_0$ and $C_{\xi}$. When $C_{\eta}\in \mathcal{C}_{\BB}$,
there are two (essentially) distinct generic possibilities.
\begin{enumerate}
\item If there is no fold catastrophe along $C_{\eta}$, then the angles $\varphi_{\eta}^{\pm}$ do not
exist and $\mathsf{W}_{\eta}$ is catastrophe-free deleted neighbourhood of the
origin in parameter space.
\item If there are $2k$ fold catastrophes along $C_{\eta}$, $k>0$, then the angles $\varphi_{\eta}^{\pm}$
are well-defined: $s_{\eta}\mp\varphi_{\eta}^{\pm}$ are the $s-$values of the
first fold catastrophes occurring respectively \textsl{before} and
\textsl{after} $\eta$ along $C_{\eta}$.
\end{enumerate}
When $C_{\eta}\in \mathcal{C}_{\infty}$, there are two (essentially) distinct
generic possibilities.
\begin{enumerate}
\item If there is no fold catastrophe along $C_{\eta}$, then the angles $\varphi_{\eta}^{\pm}$ are well
defined and $s_{\eta}\mp\varphi_{\eta}^{\pm}$ are the $s-$values of the
$\infty-$catastrophes occurring respectively \textsl{before} and \textsl{after}
$\eta$ \textsl{via} $C_{\eta}$.
\begin{figure}[t]
\begin{center}
\includegraphics[width=350pt]{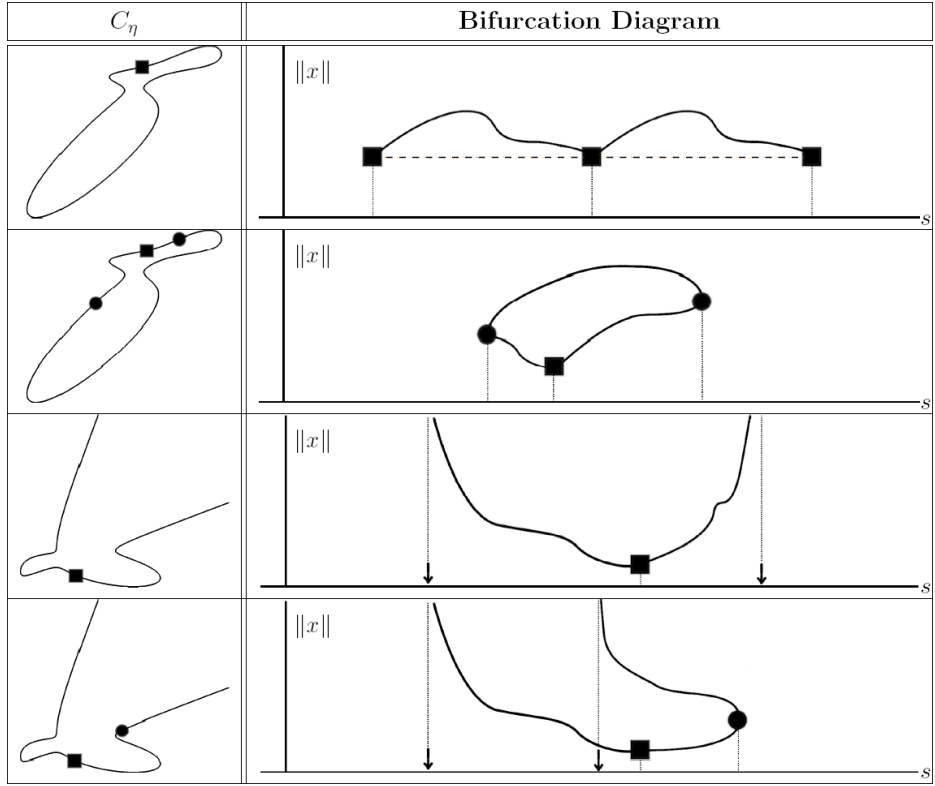}
\caption{Partial bifurcation diagrams of $P_{\rho}$ when $C_{0}\neq C_{\xi}$.
Only one branch is shown. The square represents the origin or $\xi$ and the
circles and arrows indicate fold and $\infty-$catastrophes,
respectively.}\end{center}\hrule
\end{figure}\item \label{catastrophe}If there are $k$ fold catastrophes along $C_{\eta}$, $k>0$,
then the angles $\varphi_{\eta}^{\pm}$ are well-defined: if all the fold
catastrophes lie on one side of $\eta$ (say $s>s_{\eta}$) along $C_{\eta}$ then
$s_{\eta}\mp\varphi_{\eta}^{\pm}$ are the $s-$values of the
$\infty-$catastrophes occurring \textsl{before} $\eta$ and the first fold
catastrophe \textsl{after} $\eta$ along $C_{\eta}$, respectively (or
\textsl{vice-versa}). Otherwise, $s_{\eta}\mp\varphi_{\eta}^{\pm}$ are the
$s-$values of the first fold catastrophes occurring respectively
\textsl{before} and \textsl{after} $\eta$ along $C_{\eta}$.
\end{enumerate}
Some corresponding qualitative bifurcation diagrams are shown in
figure~\thefigure.\normalsize\subsubsection{The case $C_0=C_{\xi}$} In this
instance, the bifurcation diagram must pass through $0$ at $s=0$ and $\xi$ at
$s=\pi/2$. When $C_{0}=C_{\xi}\in \mathcal{C}_{\BB}$, the number of fold
catastrophes along the curve is even; there are then three (essentially)
distinct generic possibilities.
\begin{enumerate}
\item If there is no fold catastrophe along $C_{0}=C_{\xi}$, then the angles $\varphi_{0}^{\pm}$ and $\varphi_{\xi}^{\pm}$
do not exist and $\mathsf{W}_{0}=\mathsf{W}_{\xi}$ are catastrophe-free deleted
neighbourhoods of the origin in parameter space.
\item If there is an odd number of fold catastrophes between the origin and $\xi$ along $C_{0}=C_{\xi}$, then the angles
$\mp\varphi_{0}^{\pm}$ and $\pi/2\mp\varphi_{\xi}^{\pm}$ are well-defined: they
are the $s-$values of the first fold catastrophes occurring respectively
\textsl{before} and \textsl{after} $0$ and $\xi$ \textsl{via} $C_{0}=C_{\xi}$.
\item If there is an even number of fold catastrophes between $0$ and $\xi$ along $C_{0}=C_{\xi}$, the situation is much as
described in (2), save for the fact that $C_0=C_{\xi}$ is not a loop in the
bifurcation diagram of $P_{\rho}$.
\end{enumerate}
When $C_{0}=C_{\xi}\in \mathcal{C}_{\infty}$, there are two (essentially)
distinct generic possibilities.
\begin{enumerate}
\item If there is no fold catastrophe along $C_{0}=C_{\xi}$, then the angles $\varphi_{0}^{\pm}$ and $\varphi_{\xi}^{\pm}$
are well defined and $s_{\eta}\mp\varphi_{\eta}^{\pm}$ are the $s-$values of
the $\infty-$catastrophes occurring respectively \textsl{before} and
\textsl{after} $0$ and $\xi$ \textsl{via} $C_{0}=C_{\xi}$.
\item If there are $k$ fold catastrophes along $C_{0}=C_{\xi}$, $k>0$, then the angles $\varphi_{0}^{\pm}$ and
$\varphi_{\xi}^{\pm}$ are well-defined: if no fold catastrophe lies between $0$
and $\xi$ along $C_{0}=C_{\xi}$ then $\mp\varphi_0^{\pm}$ and
$\pi/2\mp\varphi_{\eta}^{\pm}$ are determined as in the case $C_{\eta}\in
\mathcal{C}_{\infty}$, item (2) (see p.~\pageref{catastrophe}). If there are
fold catastrophes between $0$ and $\xi$ along $C_0=C_{\xi}$, then $\varphi_0^+$
and $\pi/2-\varphi_{\xi}^-$ are the $s-$values of the first fold catastrophes
occurring respectively \textsl{after} $0$ and \textsl{before} $\xi$ along
$C_{\xi}$.
\end{enumerate}
Some corresponding qualitative bifurcation diagrams are shown in
figure~\addtocounter{figure}{1}\thefigure\addtocounter{figure}{-1}.
\normalsize\subsection{The General Mapping}\label{GM}
\normalsize\sectionmark{The General Mapping} The mapping (\ref{themapping}) is
not the most general mapping satisfying (P1)$-$(P3); one should instead study
maps of the form
\begin{equation}\label{themapping2}
\mathcal{P}(x,\lambda)=x+2\pi\big[\lambda_1\mathcal{F}_0(x,\lambda_1)+\lambda_1\lambda_2\mathcal{J}(x,\lambda)
+\lambda_2\mathcal{G}_{\xi}(x,\lambda_2)\big],
\end{equation} \begin{figure}[t]\label{biffdigg}
\begin{center}
\includegraphics[width=350pt]{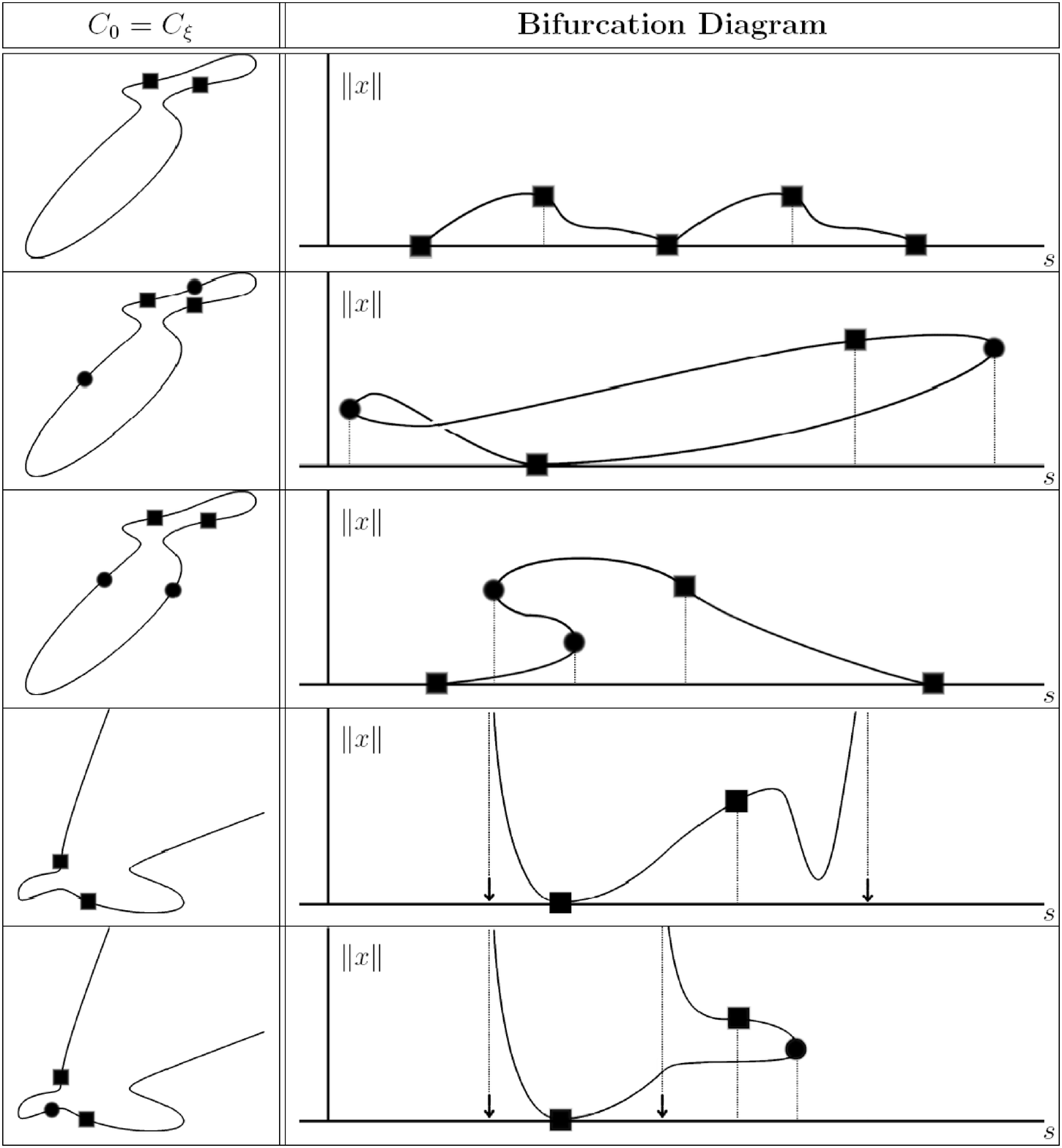}
\caption{Partial bifurcation diagrams of $P_{\rho}$ when $C_0=C_{\xi}$. The
squares represent $0$ and $\xi$, and the circles and arrows indicate fold and
$\infty-$catastrophes, respectively. The apparent self-intersection is an
artifact of the projection on the $\|x\|-s$ plane; it does not, in fact,
occur.}\end{center}\hrule
\end{figure} where $\xi\neq 0\in \GR^2$, $\mathcal{F}_0$, $\mathcal{J}$ and $\mathcal{G}_{\xi}$
are real analytic in their variables and the jacobians
$D_x\mathcal{F}_0(0,\lambda_1)$ and
 $D_x\mathcal{G}_{\xi}(\xi,\lambda_2)$ have the particular form prescribed by proposition~\ref{themappingprop2},
 which is analogous to proposition~\ref{themappingprop}.
\begin{prop}\label{themappingprop2}
If $\mathcal{F}(0,\lambda_1)\equiv 0$, $\mathcal{G}_{\xi}(\xi,\lambda_2)\equiv
0$, and if
$$D_x\mathcal{F}_0(0,\lambda_1)=\begin{pmatrix}a(\lambda_1) & -b(\lambda_1)
\\ b(\lambda_1) & a(\lambda_1)\end{pmatrix}\quad\mbox{and}\quad
D_x\mathcal{G}_{\xi}(\xi,\lambda_2)=\begin{pmatrix}c(\lambda_2) & -d(\lambda_2)
\\ d(\lambda_2) & c(\lambda_2)\end{pmatrix},$$ where $a,b,c,d:\GR\to \GR$ are
continuous in their variables and $a(0),c(0)\neq 0$, then there exists
$\omega_*>0$ such that the map defined by $(\ref{themapping2})$ satisfies
conditions \textsc{(P1)$-$(P3)}.\end{prop} \noindent Define
$\mathcal{A}:\GR^2\times \GR^2\to \GM_2(\GR)$ by
\begin{equation}\label{themA}\mathcal{A}(x,\lambda)=\begin{bmatrix} \mathcal{F}_0(x,\lambda_1)+\frac{\lambda_2}{2}\mathcal{J}(x,\lambda) &
 \frac{\lambda_1}{2}\mathcal{J}(x,\lambda)+\mathcal{G}_{\xi}(x,\lambda_2)\end{bmatrix}; \end{equation} fixed
 points of (\ref{themapping2}) are then in one-to-one correspondence with solutions of
 \begin{equation}\mathcal{A}(x,\lambda)\cdot\lambda=0. \label{whatisthisequation}\end{equation}
 Set $\mathfrak{F}(x,\lambda)=\det \mathcal{A}(x,\lambda)$. According to Taylor's theorem, there are appropriate functions
 $K_{10}$, $K_{01}$ such that $\mathfrak{F}(x,\lambda)=
 \mathfrak{F}(x,0)+\lambda_1K_{10}(x,\lambda)+\lambda_2K_{01}(x,\lambda).$
 \par \normalsize Let $\hat{x}$ be such that $\mathfrak{F}(\hat{x},0)=0$, $\det
D_x\mathfrak{F}(\hat{x},0)\neq 0$ and $\mathcal{A}(\hat{x},0)\neq 0$. Then, by
the implicit function theorem, there is a neighbourhood $\mathfrak{V}\subseteq
\GR^2$ of the origin and a unique analytic function
$\mathfrak{X}:\mathfrak{V}\to\GR^2$ such that $\mathfrak{X}(0)=\hat{x}$,
$\mathfrak{F}(\mathfrak{X}(\lambda),\lambda)\equiv 0$ and $\rank
\mathcal{A}(\mathfrak{X}(\lambda),\lambda)=1$ $\forall\, \lambda\in
\mathfrak{V}$. \newpage\noindent Define $\mathfrak{L}_{\hat{x}}=\{\lambda\in
\mathfrak{V}:\mathcal{A}(\mathfrak{X}(\lambda),\lambda)\cdot\lambda=0\}.$ A
simple rank argument shows that $\mathfrak{L}_{\hat{x}}$ is defined
\textsl{via} a single equation in two real variables, with a regular solution
at the origin\label{pregularpoint}; consequently, as a manifold,
$\mathfrak{L}_{\hat{x}}$ is one-dimensional.
\par Let $\mathcal{L}_{\hat{x}}=\ker \mathcal{A}(\hat{x},0)$. Then, there is a
small neighbourhood $\mathfrak{U}\subseteq B(0,\omega_*)$ of the origin in
parameter space for which
$$\{(\mathfrak{X}(\lambda),\lambda):\lambda\in \mathfrak{U}\cap
\mathfrak{L}_{\hat{x}}\}\quad\mbox{is a deformation of}\quad
\{(\hat{x},\lambda):\lambda\in \mathfrak{U}\cap \mathcal{L}_{\hat{x}}\}:$$ both
`curves' can  be parameterized by the same $\lambda_j$, $j=1,2$. \noindent The
preceding discussion shows that the fixed points of (\ref{themapping2}) are in
one-to-one correspondence with the fixed points of the (already studied)
truncated map
\begin{align}\label{whatsthisone}
\mathcal{P}_{\T}(x,\lambda)=x+2\pi\big[\lambda_1\mathcal{F}_0(x,0)+\lambda_2\mathcal{G}_{\xi}(x,0)\big].
\end{align}
Fold bifurcations persist under small perturbations \cite{GH,WI}. Similarly, a
generic unbounded curve remains unbounded under small
perturbations\label{pconicrazy}.\par Indeed, in the real projective plane, an
element of $\mathcal{C}_{\infty}$ meets the line at infinity in two points.
Generically, these two points are distinct and a small perturbation will not
change that fact, \textsl{i.e} the perturbed curve is still an element of
$\mathcal{C}_{\infty}$. In the non-generic case where the two points at
infinity are equal, a small perturbation will either cause the points to
separate or to vanish entirely (reminescent of a fold bifurcation of points at
infinity), \textsl{i.e} the perturbed curve either stays in
$\mathcal{C}_{\infty}$ or becomes finite.\par Thus, catastrophes generically
persist: as a result, the bifurcation diagrams of (\ref{themapping2}) and
(\ref{whatsthisone}) are (locally) topologically equivalent for small parameter
values $\lambda$. Consequently, (\ref{themapping2}) has wedge-like regions
$\mathfrak{W}_{\eta}$ corresponding to the wedge regions $\mathsf{W}_{\eta}$ of
(\ref{themapping}).
\par Finally, note that since $P$ is the `linearization' of
$\mathcal{P}$ at the origin with respect to~$\lambda$, the wedge regions
$\mathsf{W}_{\eta}$ of (\ref{themapping}) provide tangential `cones' for the
corresponding wedge-like regions $\mathfrak{W}_{\eta}$ of (\ref{themapping2}).


\normalsize\section{Numerical Simulations and Examples}

In this section, we illustrate and interpret the results of the preceding
sections through various examples. As such, the emphasis lies with qualitative
observations rather than with precise numerical analysis. First, we study
systems of PDE from a (naive) numerical perspective: we observe spiral
anchoring, as well as hysteresis and homotopy of the spiral tip. Finally, we
provide a few examples of mappings of the form (\ref{themapping2}) together
with their zero-level sets and partial bifurcation diagrams.
\normalsize\subsection{PDE, FESB and Semi-Flows}\label{PDEFESBSF} In this
section, we examine systems of partial differential equations giving rise to
semi-flows satisfying the FESB equivariance described in section 2.
\newpage\noindent The computations are carried out on a two-dimensional square
domain $[-30,30]^2$ with 200 grid points to a side and time-step $\Delta
t=0.005$ and Neumann boundary condition, using a 5-point Laplacian and i) an
explicit Runge-Kutta $2-$stage method of order two in section~\ref{SA}, and ii)
Matsui's fourth-order Runge-Kutta code based on Barkley's \texttt{EZ-Spiral} in
section~\ref{HaH}. Throughout, centers of anchoring are found \textsl{via} fast
Fourier transforms of the tip data. \normalsize\subsubsection{Spiral
anchoring}\label{SA} Consider the following small perturbation of the
FitzHugh-Nagumo equations:
\begin{figure}[t]\begin{center}
\includegraphics[height=210pt]{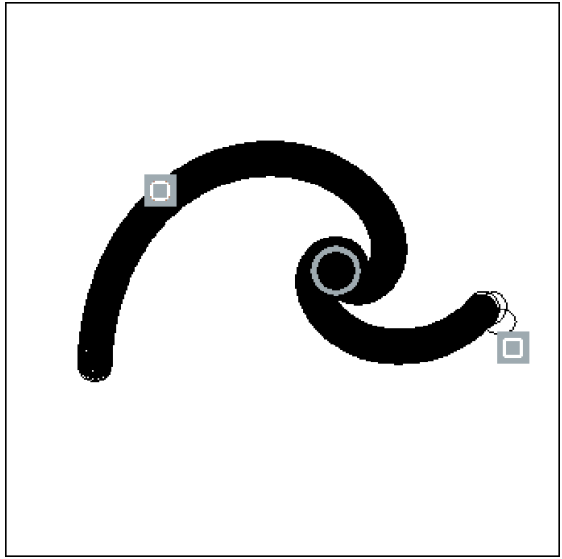}
\caption{Anchoring in (\ref{FHNM2}) with perturbations as in (\ref{thephi2}).
The spiral tip paths are plotted in black, the anchored perturbed rotating wave
is shown in gray, and the squares indicate the location of the perturbation
centers.}
\end{center}\hrule\end{figure}\begin{align}\label{FHNM2}
\begin{split}
  u_t&= \textstyle{\frac{1}{\varsigma}
  \left(u-\frac{1}{3}u^3-v\right)+\phi_1+\Delta u},\\
  v_t&=\varsigma(u+\beta-\gamma v-\phi_2),
\end{split}
\end{align} where
\begin{align}\label{thephi2}\phi_j(x) =\sqrt{2}\cos(0.05\pi) 0.12 f(x_1-c_{1,j},x_2-c_{2,j}),\ j=1,2,\end{align} $c_{1,1}=9$, $c_{2,1}=0$, $c_{1,2}=-10$, $c_{2,2}=5\sqrt{3}$ and
$$f(x)=\exp\left(-0.00086\left(x_1^2+x_2^2\right)\right).$$ Each $g_j(x)$, alone, breaks translational symmetry but preserves
rotational symmetry about $(9,0)$ (for $j=1$) or $(-10,5\sqrt{3})$ (for $j=2$).
Note that both perturbations are uniformly bounded on $\GR^2$ and that they go
to $0$ as $\|x\|\to \infty$. Under these conditions, the flow of (\ref{FHNM2})
near a normally hyperbolic rotating wave is equivalent to the flow of some
center bundle equation (\ref{system1}). Thus, if spiral waves anchor at all,
they will generically do so away from either perturbation center. This is
confirmed in figure~\thefigure, in which the transients anchor at what would be
an otherwise unremarkable location. \newl We now present the results of
simulations on a reaction-diffusion system with 4 TSB perturbations. Set
\begin{align}\label{FHNM4}
\begin{split}
  u_t&= \textstyle{\frac{1}{\varsigma}
  \left(u-\frac{1}{3}u^3-v\right)+\phi_1+\Delta u},\\
  v_t&=\varsigma(u+\beta-\gamma v+\phi_2),
\end{split}
\end{align}
where $\varsigma=0.3$, $\beta=0.6$, $\gamma=0.5$,  and where \( \phi_1 \),\(
\phi_2 \) are inhomogeneous terms which depend on $x\in \GR^2$ and are defined
by
\begin{align*}
  \phi_1(x) & = g_1(x)+g_2(x)=0.12 f_1(x_1-9,x_2)-0.10f_2(x_1+1,x_2-10), \\
  \phi_2(x) & = g_3(x)+g_4(x)=-0.12 f_1(x_1+10,x_2-5\sqrt{3})+0.08 f_3(x_1-10,x_2-10),
\end{align*}
where $A_1=0.12$, $A_2=-0.10$, $B_1=-0.12$, $B_2=0.08$,
$$f_{j}(x)=\exp\left(a_j(x_1^2+x_2^2)\right),\quad j=1,2,3,$$ $a_1=-0.00086$, $a_2=-0.0008$ and $a_3=-0.0009$.\par Each $g_j(x)$, alone,
breaks translational symmetry but preserves rotational symmetry about
$c_1=(9,0)$ (for $j=1$), $c_2=(-1,10)$ (for $j=2$), $c_3=(-10,5\sqrt{3})$ (for
$j=3$) and $c_4=(10,10)$ (for $j=4$). Note that the four perturbations are
uniformly bounded on $\GR^2$ and that they go to $0$ as $\|x\|\to \infty$. As
predicted, anchoring takes place away from the $c_j$,
\begin{figure}[t]
\begin{center}
\includegraphics[width=315pt]{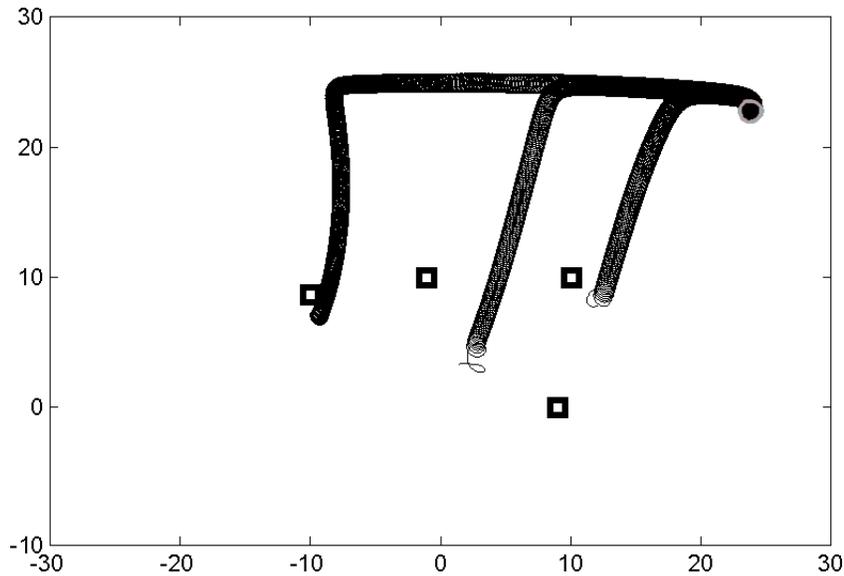}
\caption{Anchoring in the FitzHugh-Nagumo equations (\ref{FHNM4}). The spiral
tip paths are plotted in black, the anchored perturbed rotating wave is shown
in gray and the black squares indicate the location of the perturbation
centers.}
\end{center}\hrule\end{figure} $j=1,\ldots,4$.\par
The transients in figure~\thefigure\ appear to first (hyperbolically) approach
some manifold along which they travel to the anchored perturbed rotating wave;
this will be the topic of an upcoming paper. \normalsize\subsubsection{Homotopy
and hysteresis of rotating waves}\label{HaH} Following \cite{MPMPV}, define the
modified Oregonator
\begin{align}\label{MOM}
\begin{split}
  u_t&= \textstyle{\frac{1}{\varsigma}
  \left(u-u^2-(fv+\phi)\frac{u-q}{u+q}\right)+\Delta u},\\
  v_t&=(u-v)+0.6\Delta v,
\end{split}
\end{align}
where \( f =1.4\), \( q=0.002 \), \( \varsigma=0.05 \) and \( \phi \) is an
inhomogeneous term which depends on $x\in \GR^2$. When \( \phi\equiv 0 \),
(\ref{MOM}) has full Euclidean symmetry. \newpage\noindent In the following
simulations, \( \phi \) is the sum of two Gaussian bells:
\begin{align*}
  \phi(x) &=\alpha_1 \exp\left(-\frac{\|(x_1,x_2)-(15,15)\|^2}{{\beta_1}^2}\right)
  +\alpha_2 \exp\left(-\frac{\|(x_1,x_2)-(18.75,15)\|^2}{{\beta_2}^2}\right),
\end{align*}
with \( \alpha_1 \), \( \alpha_2 \), \( \beta_1 \), \( \beta_2 \in \mathbb{R}
\), and $\beta_1,\beta_2\neq 0$. Each $g_j(x)$, alone, breaks translational
symmetry but preserves rotational symmetry about $c_1=(15,15)$ (for $j=1$) or
$c_2=(18.75,15)$ (for $j=2$). Note that both perturbations are uniformly
bounded on $\GR^2$ and that they go to $0$ as $\|x\|\to \infty$. \par When
$\alpha=(\alpha_1,0)\neq 0$, (\ref{MOM}) is $\GSO(2)_{c_1}\!-$equi\-va\-riant.
Similarly, (\ref{MOM}) is $\GSO(2)_{c_2}\!-$ equi\-va\-riant when
$\alpha=(0,\alpha_2)\neq 0$ and trivially equivariant when
$\alpha_1,\alpha_2\neq 0$. \newl Set $\beta_1=\beta_2=1$ and $\rho_*=0.01$.
Along the path $\alpha(\tau)=\gamma_1(\tau)=\rho_*\left( \cos(\tau),\sin(\tau)
\right)^{\!\top}$ in parameter space, (\ref{MOM}) undergoes a homotopy of
perturbed rotating waves, whose tip paths deform continuously from a circle
centered at $c_1$, when $\tau=0$, to a circle centered at $c_2$, when
$\tau=\pi/2$
\begin{figure}[t]
\begin{center}
\includegraphics[width=250pt]{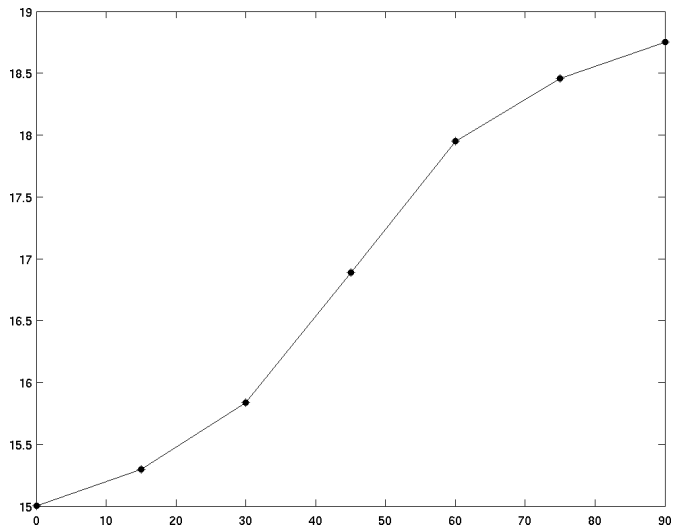}
\caption{Homotopy of the spiral tip path in (\ref{MOM}). The first spatial
coordinates of the anchoring centers are plotted against $\tau$; compare this
image with the first bifurcation diagram on p.~\pageref{biffdigg}.}
\end{center}\hrule\end{figure} (see figure~\thefigure).
\newl Along the path $\alpha(\tau)=\gamma_2(\tau)=\frac{1}{10}\gamma_1(\tau)$, however,
the homotopy is replaced by hysteresis. As the parameters vary along the path,
(\ref{MOM}) has an anchored perturbed rotating wave whose tip path deforms
continuously from a circle centered at $c_1$. At $\aleph_1\in \gamma_2$, the
rotating wave jumps (discontinuously) to another anchored perturbed rotating
wave, whose tip path deforms continuously from a circle centered at $c_2$.
\newpage\noindent Following $\gamma_2$ in the opposite direction leads to similar
behaviour, this time with the discontinuous jump taking place at
$\aleph_2\in\gamma_2$, as can be seen in
 \addtocounter{figure}{1}figure~\thefigure\addtocounter{figure}{-1}. Consequently, there must be a third unstable rotating wave
(which escapes detection by direct means) appearing and then disappearing at
$\aleph_1$ and $\aleph_2$, respectively, in saddle-node bifurcations.\newl In a
sense, both of these occurrences have been predicted by the analysis provided
in section~\ref{CSA}; consult, for instance, the first and third bifurcation
diagrams on p.~\pageref{biffdigg}. \normalsize\subsection{Wedges and
Catastrophes}\label{tcn2} \begin{figure}[t]
\begin{center}
\includegraphics[width=250pt]{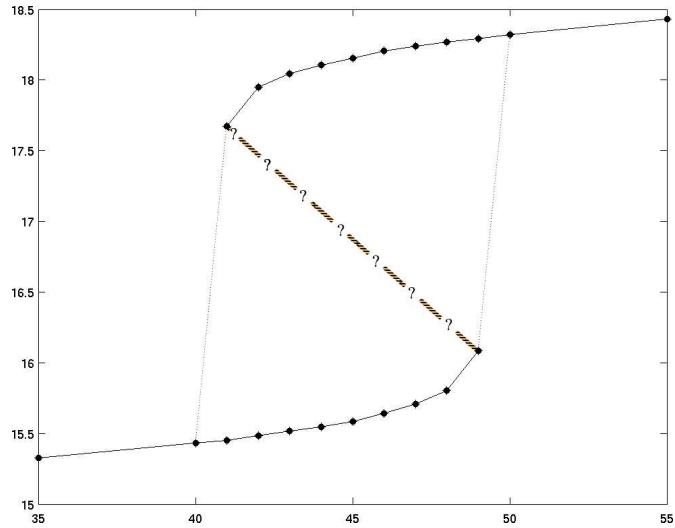}
\caption{Hysteresis of the spiral tip path in (\ref{MOM}). The first spatial
coordinates of the anchoring centers are plotted against $\tau$. The question
marks interpolate (roughly) the unstable rotating waves. Compare this image
with the third bifurcation diagram on p.~\pageref{biffdigg}.}
\end{center}\hrule\end{figure}In this final section, we provide a partial catalogue
of (partial) bifurcation diagrams for mappings of the form
\begin{align}\label{thePP}P(x,\lambda)=x+2\pi\big[\lambda_1F_0(x)+\lambda_2G_{\xi}(x)\big],\end{align}
where $\lambda\in \GR^2$, $\xi =\left(2, 2\right)^{\!\top}$,
\begin{align*}F_0(x)&=\begin{pmatrix}2x_1-x_2+\sum a_{i,j}x_1^ix_2^j \\
x_1+2x_2+\sum b_{i,j}x_1^ix_2^j\end{pmatrix}
\cdot f_0(x),\\
G_{\xi}(x)&=\begin{pmatrix}7-3x_1-\frac{x_2}{2}+\sum c_{i,j}(x_1-2)^i(x_2-2)^j
\\
5-3x_1+\frac{x_2}{2}+\sum d_{i,j}(x_1-2)^i(x_2-2)^j\end{pmatrix}\cdot
g_{\xi}(x),\end{align*} $a_{i,j},b_{i,j},c_{i,j},d_{i,j}\in \GR$, $i+j>1$, and
$f_0$ and $g_{\xi}$ are \normalsize continuous functions such that
(\ref{thePP}) satisfies proposition \ref{themappingprop}; we can then use the
visual criterion of section~\ref{TVC} to understand the nature of the
bifurcation diagram of the associated map $P_{0.01}$. Furthermore, the wedge
angles can be read directly from the bifurcation diagram. \newl In the figures
of this section, $C_0$ and $C_{\xi}$ are shown in black or gray. Fold
catastrophes are indicated by circles, $\infty-$catastrophes by arrows and the
squares mark both the origin and $\xi$. \normalsize\subsubsection{The
Elowyn-Bonhomme map} \label{EBmappp} The Elowyn-Bonhomme (EB) map is obtained
from (\ref{thePP}) by setting $f_0\equiv g_{\xi}\equiv 1$,
$a_{1,1}=b_{2,0}=c_{0,2}=d_{1,1}=d_{2,0}=1$, $a_{0,2}=-1$ and all other
coefficients to $0$; the corresponding $\kappa(\mathfrak{Z})$ is shown in
figure~\thefigure. In this instance, $C_0\in\mathcal{C}_{\BB}$ and
$C_{\xi}\in\mathcal{C}_{\infty}$. For the EB map, $\omega_*=\frac{12}{37\pi}$.
Let $\rho=0.01<\omega_*$. Using a pseudo-arc length continuation algorithm (see
\cite{KE} for details), a partial bifurcation diagram of $P_{\rho}$ (ignoring
all fixed point branches but those through $\eta$ at $s=s_{\eta}$, for $\eta\in
\{0,\xi\}$) is built: the results can be seen in
 figure\addtocounter{figure}{1}~\thefigure.\addtocounter{figure}{-1}\par There are 6 fold catastrophes: two along $C_0$ and four
along $C_{\xi}$. Their location can be recovered
 directly from $\kappa(\mathfrak{Z})$
and $\mathcal{R}_j$,
 $j=1,2$, with the help of proposition
\ref{therealprop}: in
figure\addtocounter{figure}{1}~\thefigure\addtocounter{figure}{-1}, the six
intersections that satisfy the appropriate hypotheses are marked with circles.
Each corresponds to one of the six fold catastro\-phes observed in
\addtocounter{figure}{-1}figure~\thefigure.\addtocounter{figure}{1}\par
Furthermore, two $\infty-$catas\-tro\-phes occur \textsl{via} $C_{\xi}$. The
interesting values for the EB map are compiled in
\begin{table}
\caption{EB map catastrophes along $C_0$ and $C_{\xi}$.}\begin{indented}\item[]
\begin{tabular}{ccccc}\br
\textit{Curve} & \textit{Type} & $x^*$ & $s^*$ & \textit{Wedge Angle} \\
 \mr
$C_0$ & Fold & $\left(1.2483,-0.1286 \right)^{\!\top}$ & $5.9809$ & $\varphi_{0}^-\approx 0.3023$\\
$C_0$ & Fold & $\left(0.2269,-3.4760 \right)^{\!\top}$ & $0.2308$ & $\varphi_{0}^+\approx 0.2308$\\
$C_{\xi}$ & Fold & $\left(0.3371,3.1473 \right)^{\!\top}$ & $1.1020$ & $\varphi_{\xi}^-\approx 0.4688$\\
$C_{\xi}$ & Fold & $\left(2.2769,0.2982 \right)^{\!\top}$ & $2.1125$ & $\varphi_{\xi}^+\approx 0.5417$\\
$C_{\xi}$ & Fold & $\left(-3.2933,6.1024 \right)^{\!\top}$ & $1.1581$ & n.a. \\
$C_{\xi}$ & Fold & $\left(5.6733,-1.2807 \right)^{\!\top}$ & $1.9267$ & n.a. \\
$C_{\xi}$ & Infinity & n.a. & $1.0172$ & n.a. \\
$C_{\xi}$ & Infinity & n.a. & $2.3562$ & n.a. \\
 \br
\end{tabular}
\end{indented}\end{table}table~\thetable\ above.~\normalsize We continue by providing examples that
highlight the various possibilities.
\begin{figure}[t]
\begin{center}
\includegraphics[width=370pt]{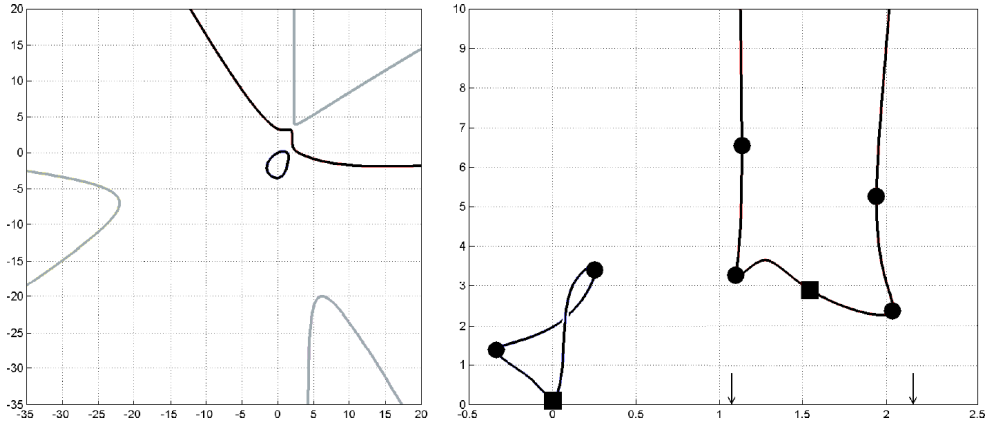}
\caption{On the left: zero-level set $\kappa(\mathfrak{Z})$ for the EB map. On
the right: fixed point branches of $P_{0.01}$ (the apparent self-intersection
on $C_0$ is due to a projection onto the $\|x\|-s$ plane). Both $C_0$ (the
curve through $0$) and $C_{\xi}$ (the curve through $\xi$) are shown in black.}
\end{center}\hrule\end{figure}
\begin{figure}[t]
\begin{center}
\includegraphics[width=200pt]{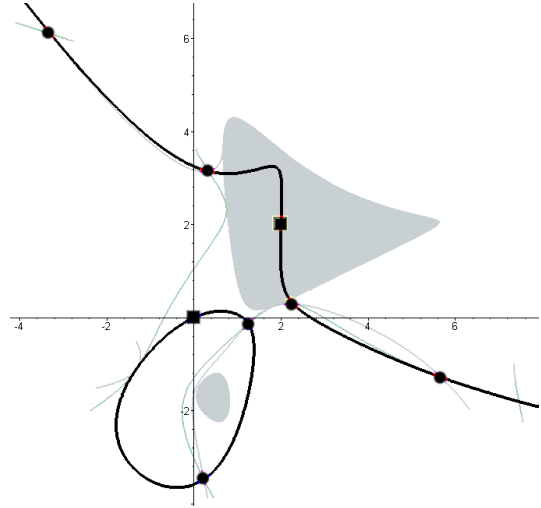}
\caption{Intersections of the zero-level sets $\kappa(\mathfrak{Z})$ (thick
black lines) and $\mathcal{R}_1$ and $\mathcal{R}_2$ (thin gray lines) of the
Elowyn-Bonhomme problem. The squares represent $0$ and $\xi$; the points that
satisfy the hypotheses of proposition~\ref{therealprop} are marked with
circles. The light gray region is part of the planar set for which
$C(x)^2-4B(x)E(x)<0$.}
\end{center}\end{figure}
\normalsize\subsubsection{The first example}
\begin{figure}[t]
\begin{center}
\includegraphics[height=142pt]{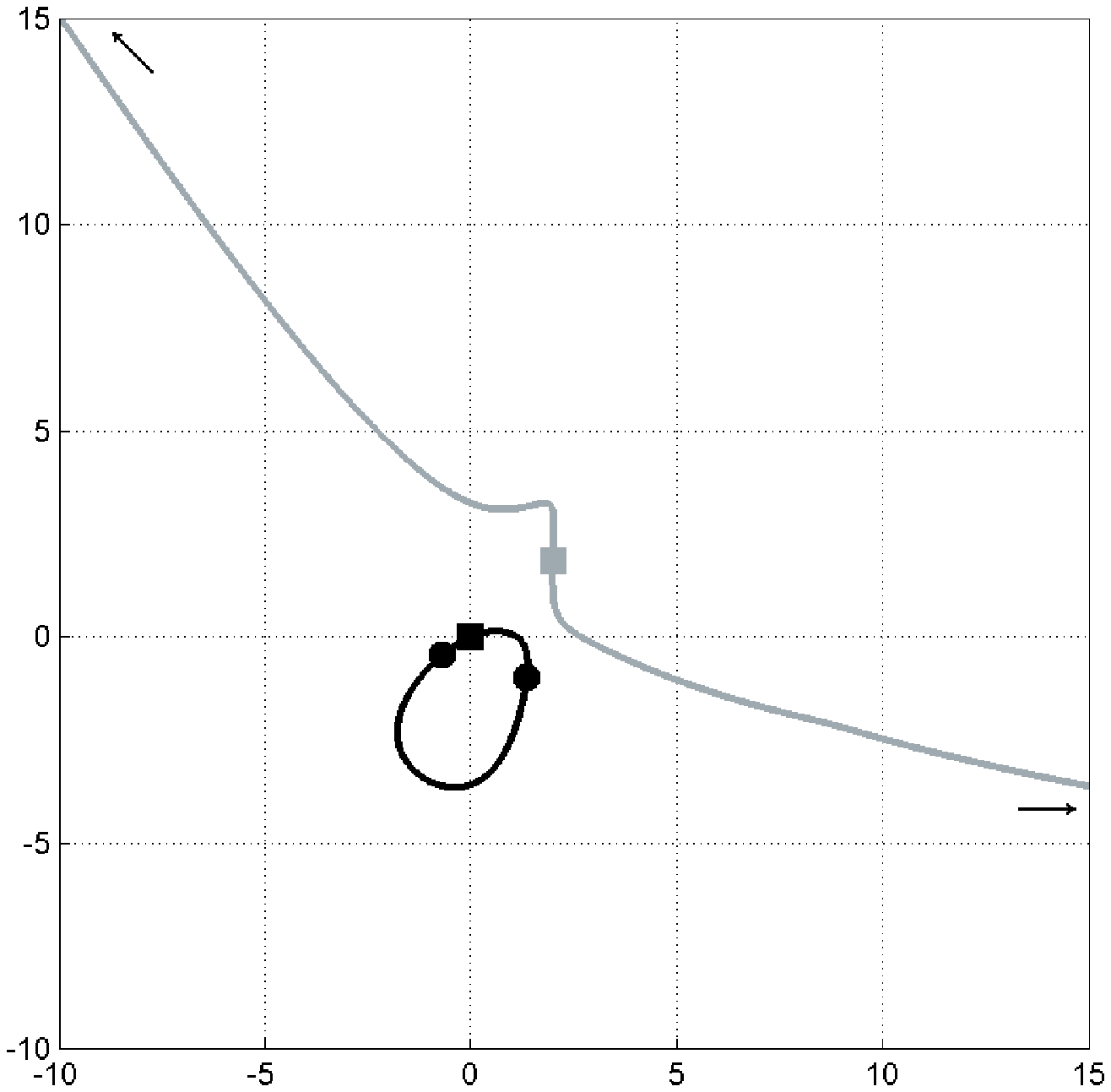}\qquad\includegraphics[height=142pt]{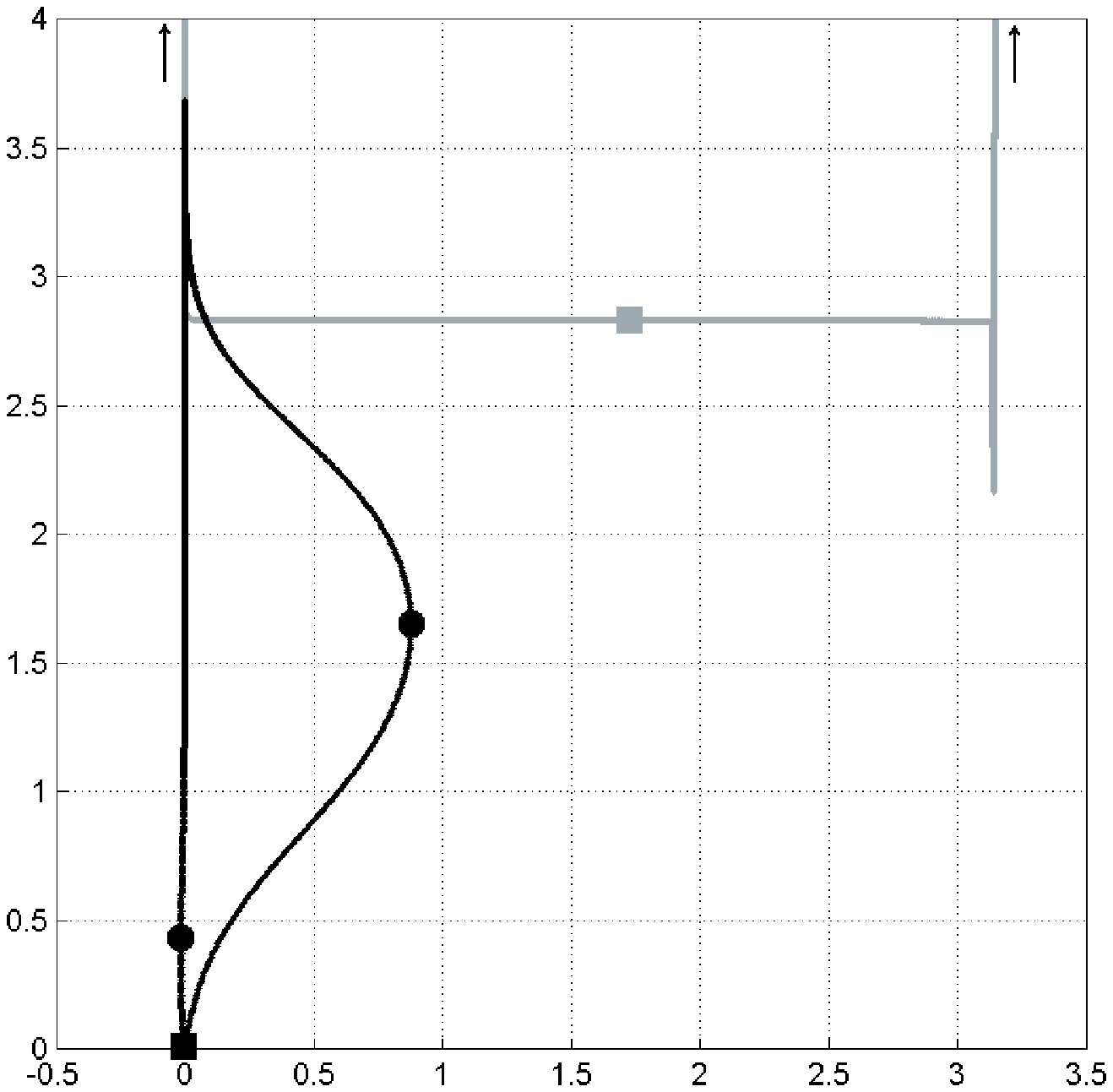}
\caption{The first example: $\kappa(\mathfrak{Z})$ (left), bifurcation diagram
of $P_{0.01}$ (right).}
\end{center}\hrule
\end{figure}
\begin{xalignat*}{6}
a_{2,0}&=0 & a_{1,1}&=1 & a_{0,2}&=-1 &b_{2,0}&= 1 &b_{1,1}&= 0 &b_{0,2}&=0 \\
c_{2,0}&=0 & c_{1,1}&=0 & c_{0,2}&=1  &d_{2,0}&= 1 &d_{1,1}&= 1 &d_{0,2}&=0
\end{xalignat*}
$$f_0(x)= \textstyle{\exp\left(-(x_1^2+x_2^2)/10\right)}\qquad g_{\xi}(x)=\textstyle{\exp\left(-\left((x_1-2)^2+(x_2-2)^2\right)/14\right)}.$$
See figure~\thefigure\ for a portion of $\kappa(\mathfrak{Z})$ and a partial
bifurcation diagram of $P_{0.01}$. By construction, the zero-level set for this
first example is exactly the zero-level set of the EB map; however, their
bifurcation diagrams are not topologically equivalent (compare with figure 9).
In this instance, the wedge angles record fold catastrophes on $C_0 \mbox{
(black)}\in \mathcal{C}_{\BB}$ and $\infty-$catastrophes on $C_{\xi}\mbox{
(gray)}\in \mathcal{C}_{\infty}$. Note further that this map provides an
instance when the anchoring wedges overlap. \normalsize\subsubsection{The
second example}
\begin{figure}[t]
\begin{center}
\includegraphics[height=142pt]{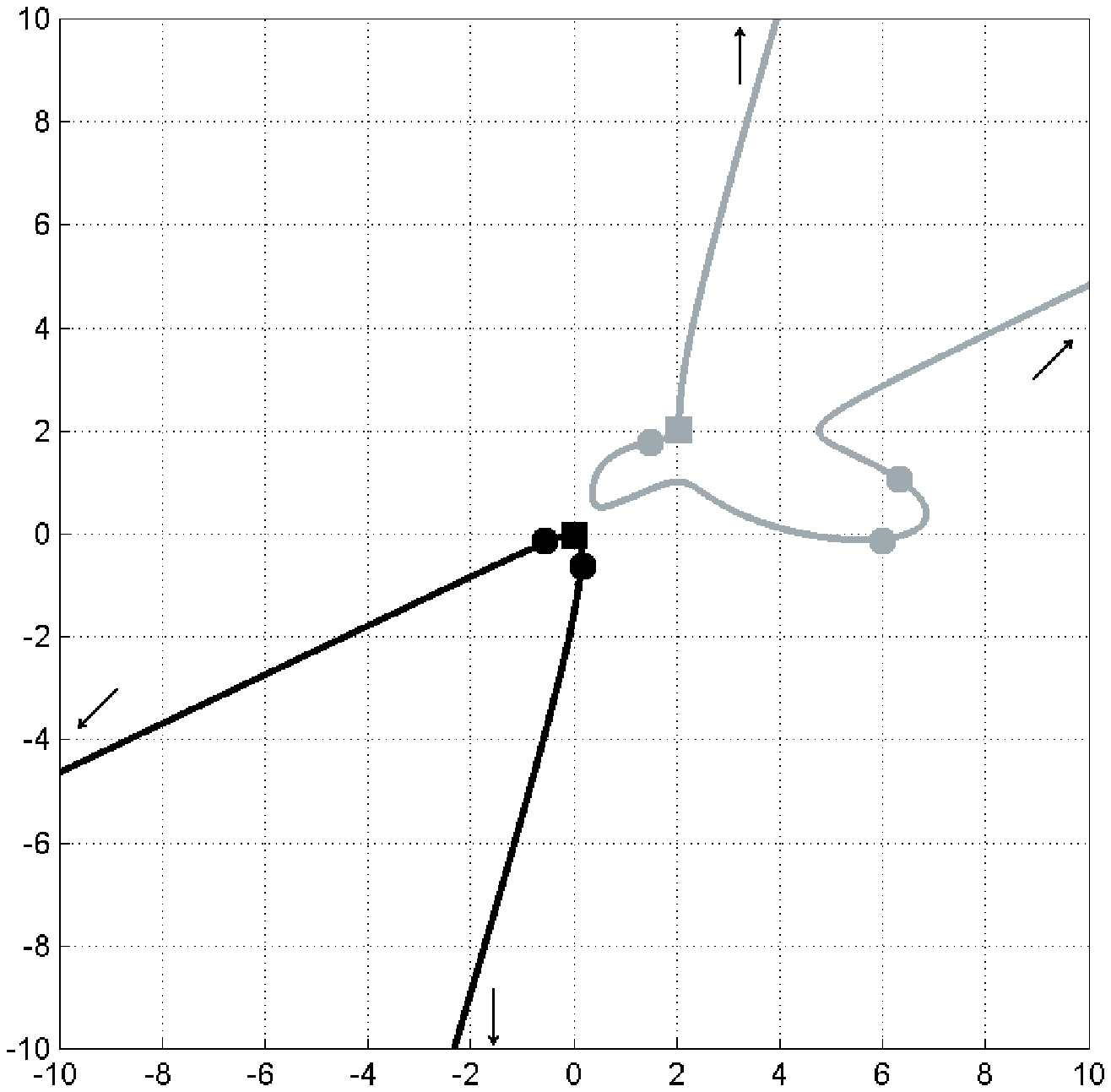}\qquad\includegraphics[height=142pt]{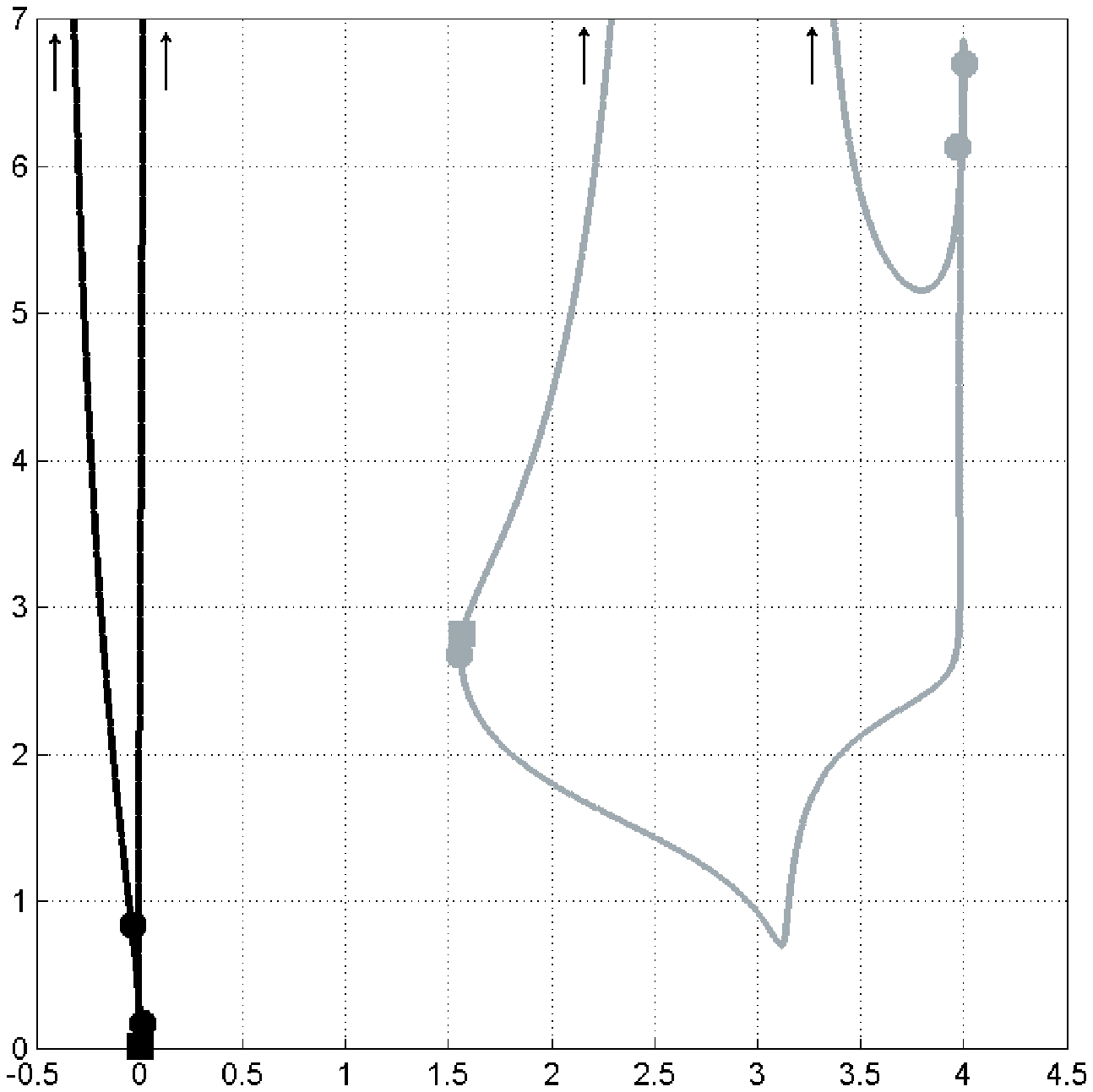}
\caption{The second example: $\kappa(\mathfrak{Z})$ (left), bifurcation diagram
of $P_{0.01}$ (right).}
\end{center}
\end{figure}
\begin{xalignat*}{6}
a_{2,0}&=\textstyle{\frac{2}{5}} & a_{1,1}&=\textstyle{\frac{43}{10}} & a_{0,2}&=-\textstyle{\frac{17}{2}} &b_{2,0}&= \textstyle{\frac{51}{10}} &b_{1,1}&= -\textstyle{\frac{91}{10}} &b_{0,2}&=-\textstyle{\frac{12}{5}} \\
c_{2,0}&=-\textstyle{\frac{16}{5}} & c_{1,1}&=-\textstyle{\frac{99}{10}} &
c_{0,2}&=-\textstyle{\frac{48}{5}}  &d_{2,0}&= -\textstyle{\frac{53}{10}}
&d_{1,1}&=\textstyle{\frac{61}{10}} &d_{0,2}&=-\textstyle{\frac{99}{10}}
\end{xalignat*}
$$f_0(x)=1\qquad g_{\xi}(x)=1$$
See figure~\thefigure\ for a portion of $\kappa(\mathfrak{Z})$ and a partial
bifurcation diagram of $P_{0.01}$. In this instance, $C_0 \mbox{
(black)},C_{\xi}\mbox{ (gray)}\in \mathcal{C}_{\infty}$, and
$\varphi_0^{\pm},\varphi_{\xi}^-$ record fold catastrophes while
$\varphi_{\xi}^+$ records an $\infty-$catastrophe.
\normalsize\subsubsection{The third example}
\begin{figure}[t]
\begin{center}
\includegraphics[height=142pt]{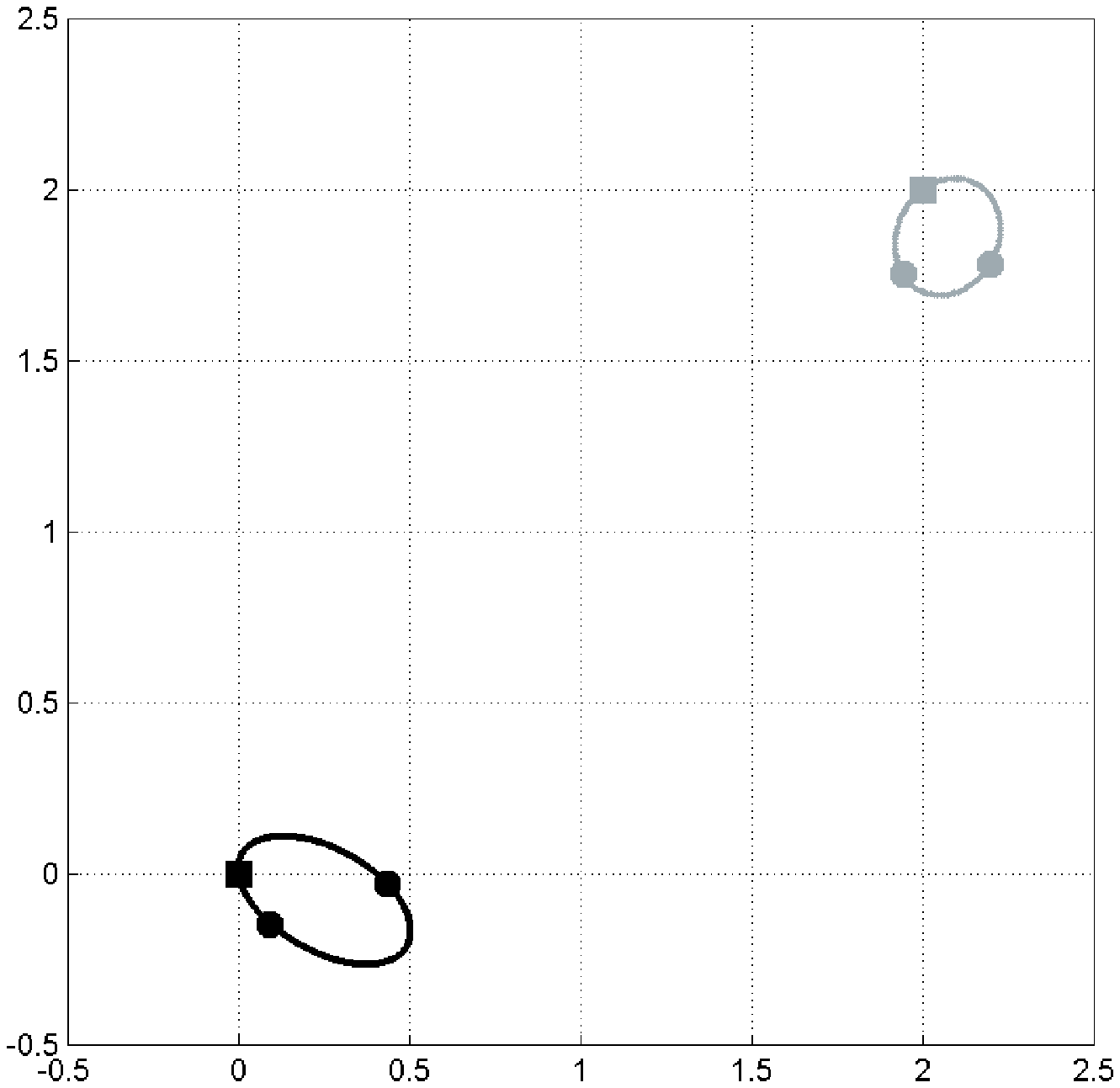}\qquad\includegraphics[height=142pt]{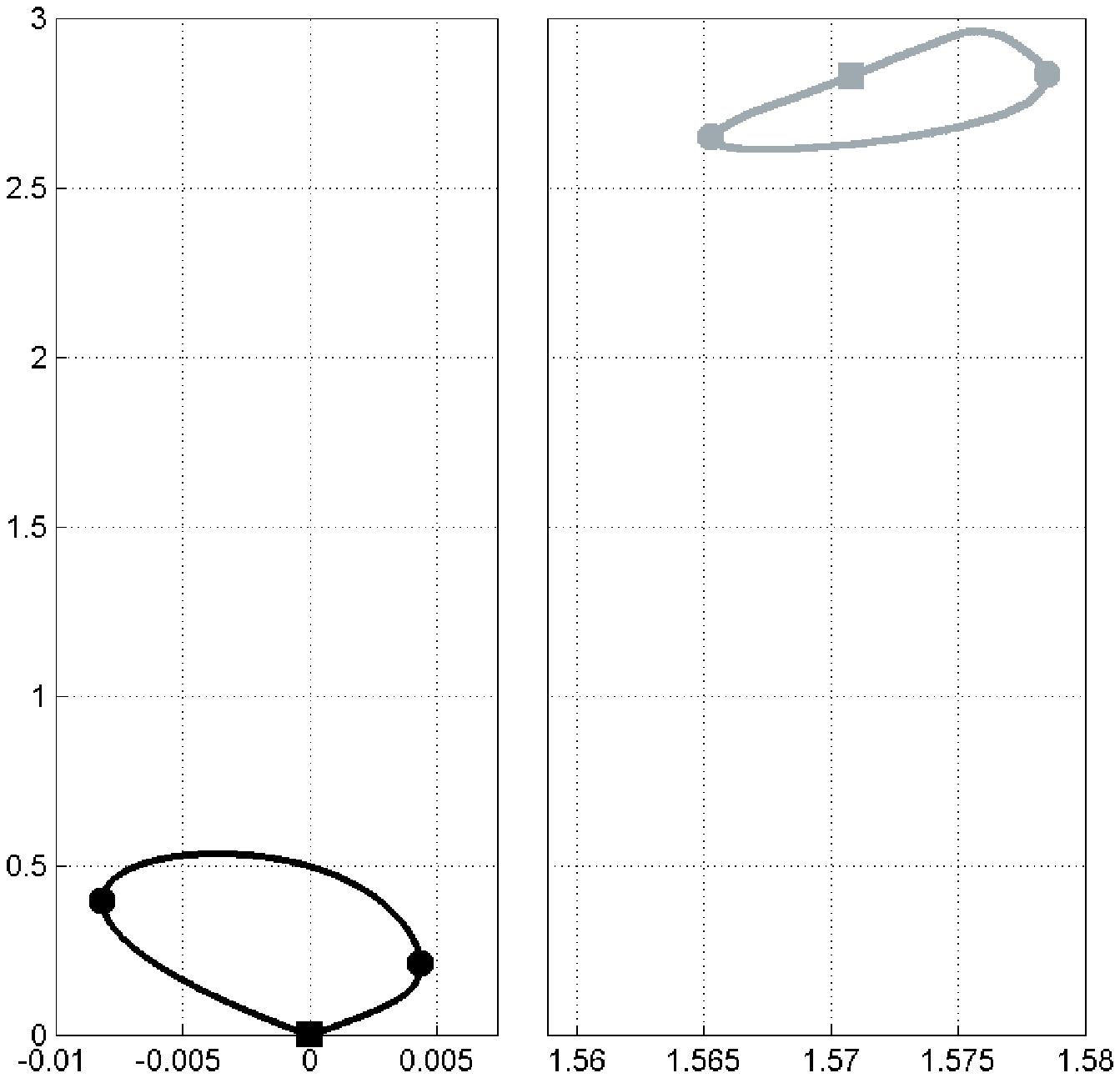} \caption{The third example: $\kappa(\mathfrak{Z})$ (left), bifurcation diagram
of $P_{0.01}$ (right).}
\end{center}\hrule
\end{figure}
\begin{xalignat*}{6}
a_{2,0}&=-\textstyle{\frac{53}{10}} & a_{1,1}&=-\textstyle{\frac{28}{5}} & a_{0,2}&=-\textstyle{9} &b_{2,0}&= -\textstyle{\frac{9}{5}} &b_{1,1}&= -\textstyle{\frac{41}{10}} &b_{0,2}&=-\textstyle{\frac{57}{10}} \\
c_{2,0}&=\textstyle{\frac{22}{5}} & c_{1,1}&=\textstyle{\frac{21}{5}} &
c_{0,2}&=-\textstyle{\frac{3}{2}}  &d_{2,0}&= -\textstyle{9}
&d_{1,1}&=\textstyle{\frac{49}{10}} &d_{0,2}&=-\textstyle{10}
\end{xalignat*}
$$f_0(x)=1\qquad g_{\xi}(x)=1$$
See figure~\thefigure\ for a portion of $\kappa(\mathfrak{Z})$ and a partial
bifurcation diagram of $P_{0.01}$. In this instance, $C_0\mbox{
(black)},C_{\xi}\mbox{ (gray)}\in \mathcal{C}_{\BB}$, and the wedge angles all
record fold catastrophes. \normalsize\subsubsection{The fourth example}
\begin{figure}[t]
\begin{center}
\includegraphics[height=142pt]{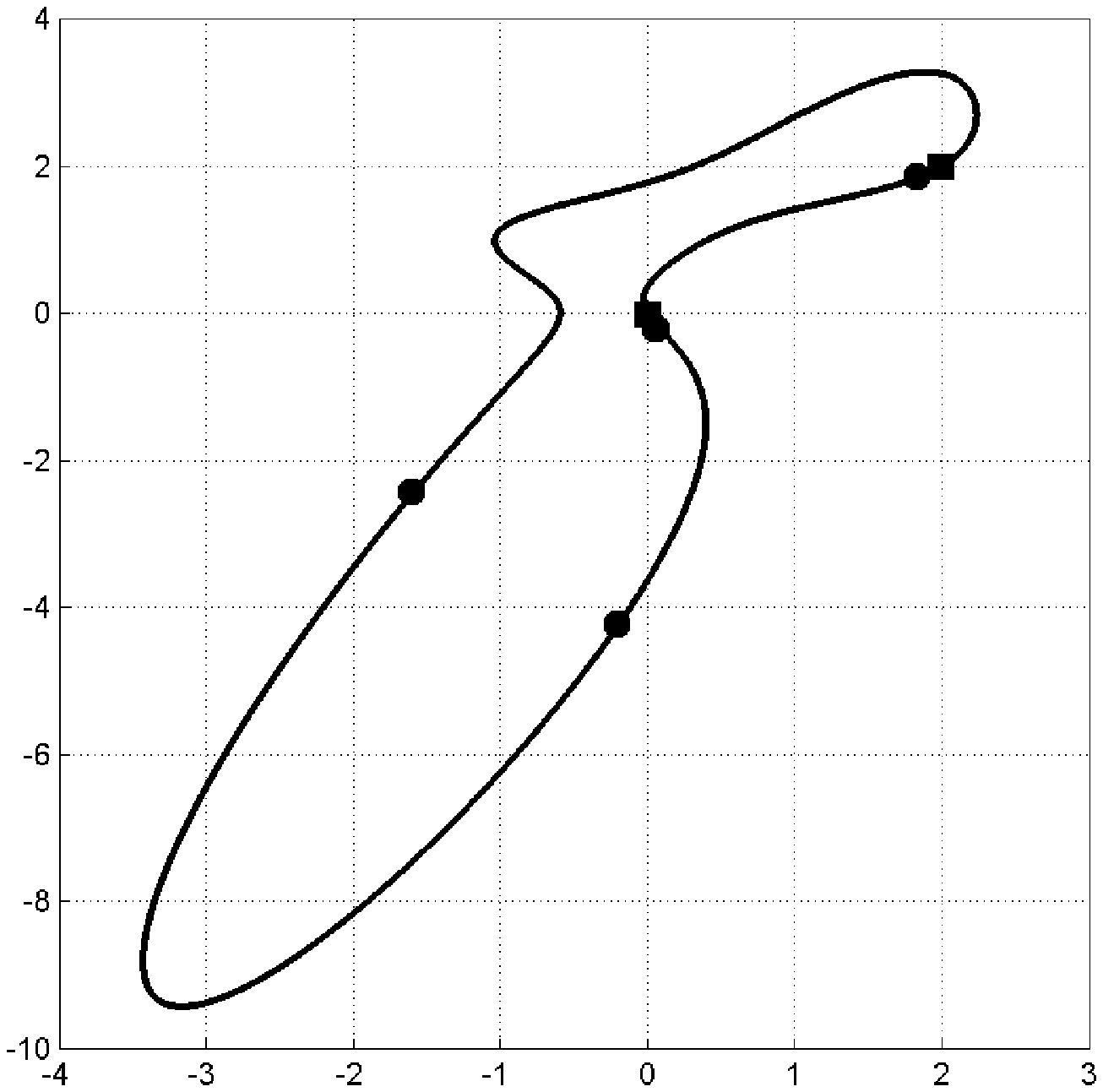}\qquad\includegraphics[height=142pt]{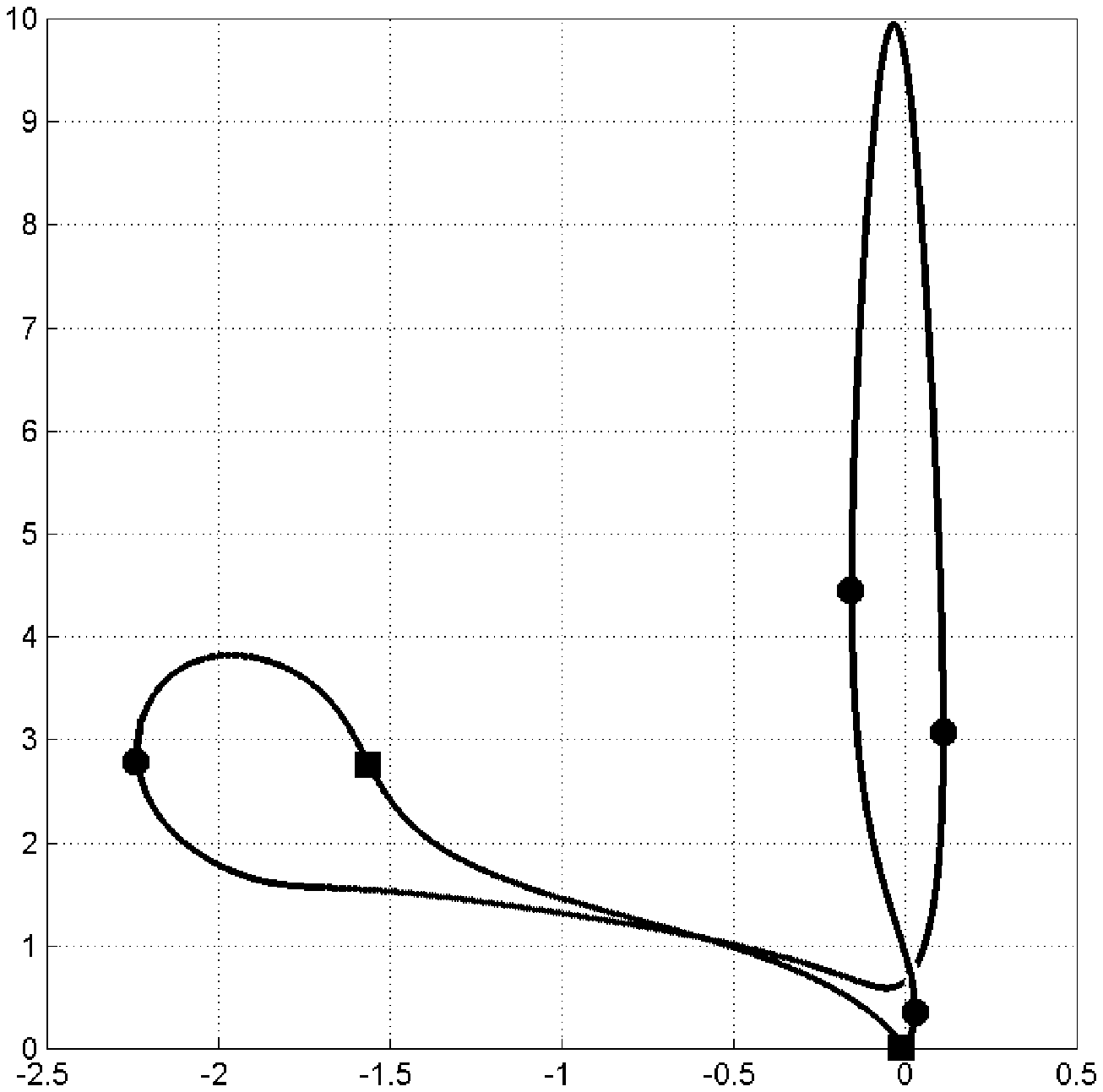} \caption{The fourth example: $\kappa(\mathfrak{Z})$ (left), bifurcation diagram
of $P_{0.01}$ (right).}\end{center}
\end{figure}
\begin{xalignat*}{6}
a_{2,0}&=\textstyle{2} & a_{1,1}&=\textstyle{\frac{13}{5}} & a_{0,2}&=-\textstyle{\frac{13}{10}} &b_{2,0}&= \textstyle{\frac{16}{5}} &b_{1,1}&= -\textstyle{\frac{21}{10}} &b_{0,2}&=\textstyle{\frac{3}{5}} \\
c_{2,0}&=-\textstyle{\frac{4}{5}} & c_{1,1}&=-\textstyle{\frac{3}{10}} &
c_{0,2}&=-\textstyle{\frac{29}{10}}  &d_{2,0}&= -\textstyle{\frac{9}{10}}
&d_{1,1}&=\textstyle{\frac{7}{10}} &d_{0,2}&=-\textstyle{\frac{3}{2}}
\end{xalignat*}
$$f_0(x)=1\qquad g_{\xi}(x)=1$$
See figure~\thefigure\ for a portion of $\kappa(\mathfrak{Z})$ and a partial
bifurcation diagram of $P_{0.01}$. In this instance, $C_0=C_{\xi}\in
\mathcal{C}_{\BB}$, the wedge angles all record fold catastrophes and the
anchoring wedges overlap. \normalsize\subsubsection{The fifth example}
\begin{figure}[t]
\begin{center}
\includegraphics[height=142pt]{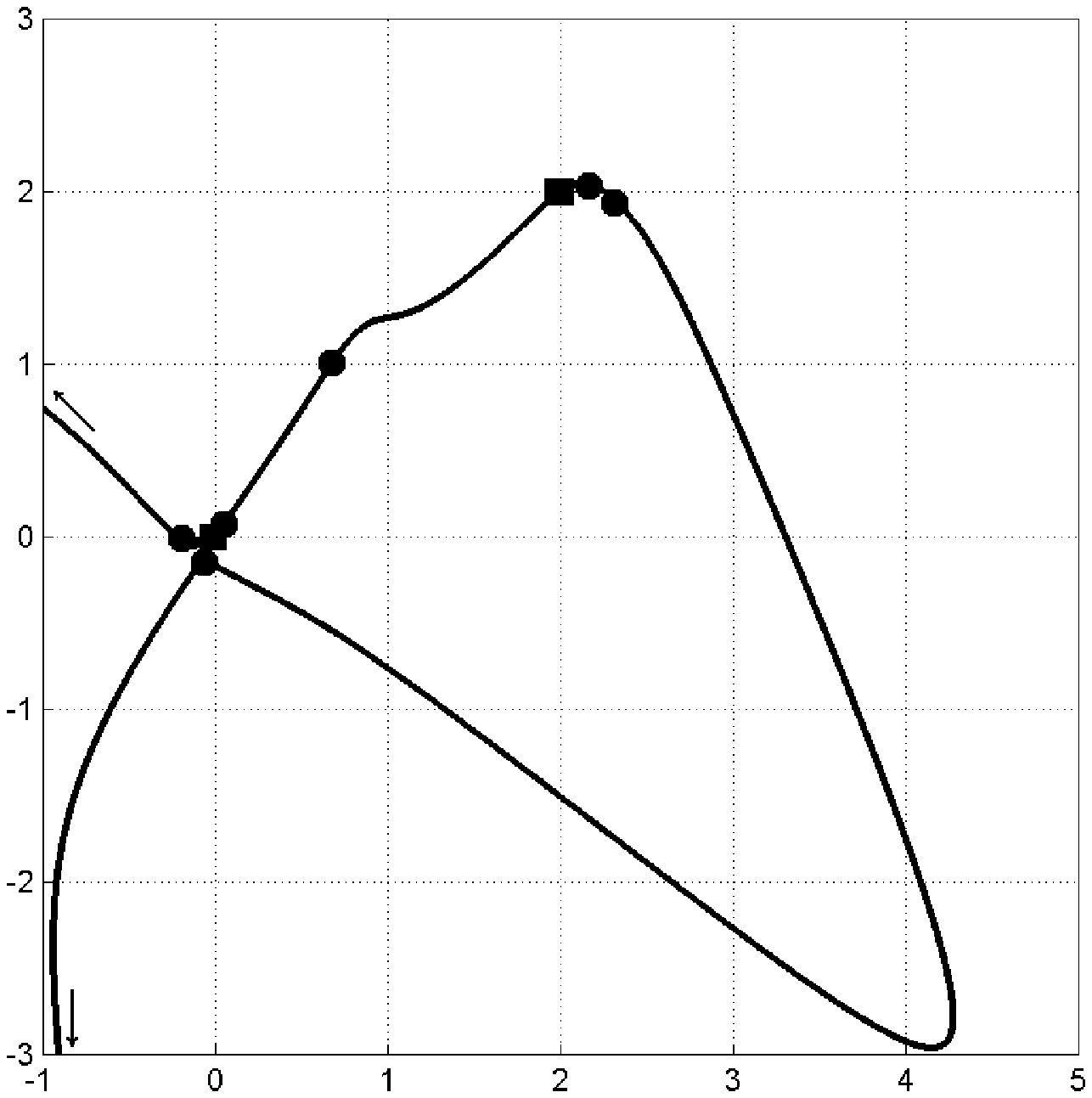}\qquad
\includegraphics[height=142pt]{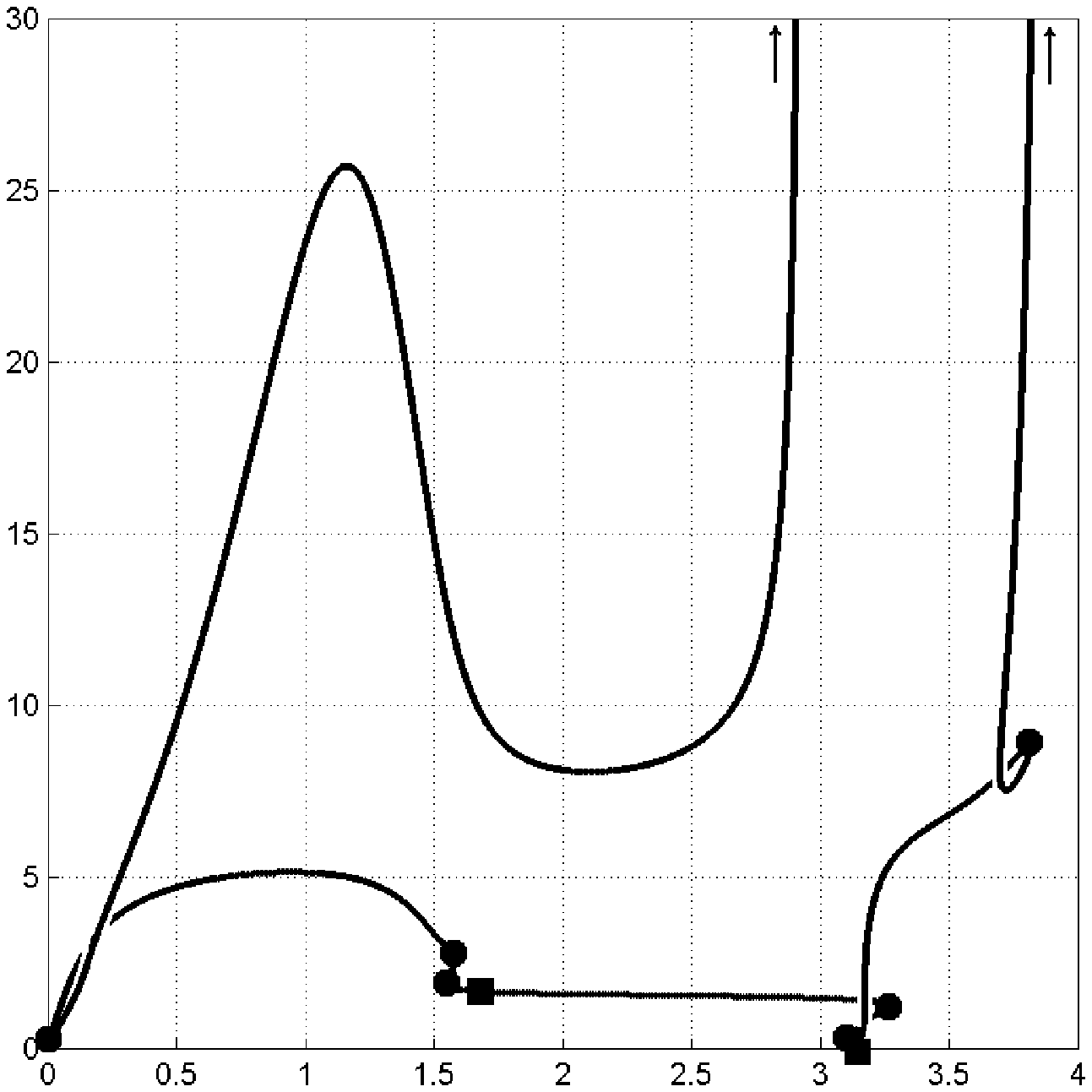} \caption{The fifth example: $\kappa(\mathfrak{Z})$ (left), bifurcation diagram
of $P_{0.01}$ (right).}
\end{center}\hrule
\end{figure}
\begin{xalignat*}{7}
a_{2,0}&=\textstyle{\frac{59}{10}} & a_{1,1}&=\textstyle{\frac{43}{10}} & a_{0,2}&=-\textstyle{7} & a_{3,0}&=\textstyle{\frac{38}{5}} & a_{2,1}&=\textstyle{\frac{27}{5}} & a_{1,2}&=-\textstyle{\frac{59}{10}} & a_{0,3}&=-\textstyle{\frac{1}{5}} \\
b_{2,0}&=-\textstyle{\frac{37}{5}} & b_{1,1}&=-\textstyle{\frac{51}{10}} & b_{0,2}&=-\textstyle{\frac{9}{5}} & b_{3,0}&=\textstyle{\frac{63}{10}} & b_{2,1}&=\textstyle{\frac{41}{10}} & b_{1,2}&=-\textstyle{4} & b_{0,3}&=\textstyle{\frac{43}{10}} \\
c_{2,0}&=\textstyle{\frac{11}{5}} & c_{1,1}&=-\textstyle{\frac{36}{5}} & c_{0,2}&=-\textstyle{\frac{38}{5}} & c_{3,0}&=-\textstyle{\frac{61}{10}} & c_{2,1}&=\textstyle{2} & c_{1,2}&=-\textstyle{\frac{4}{5}} & c_{0,3}&=-\textstyle{\frac{8}{5}} \\
d_{2,0}&=-\textstyle{\frac{49}{10}} & d_{1,1}&=-\textstyle{\frac{38}{5}} &
d_{0,2}&=\textstyle{\frac{27}{10}} & d_{3,0}&=-\textstyle{\frac{79}{10}} &
d_{2,1}&=\textstyle{\frac{69}{10}} & d_{1,2}&=-\textstyle{\frac{29}{5}} &
d_{0,3}&=-\textstyle{\frac{41}{10}}
\end{xalignat*}
$$f_0(x)=1\qquad g_{\xi}(x)=1$$
See figure~\thefigure\ for a portion of $\kappa(\mathfrak{Z})$ and a partial
bifurcation diagram of $P_{0.01}$. In this instance, $C_0=C_{\xi}\in
\mathcal{C}_{\infty}$ and the wedge angles all record fold catastrophes.
\normalsize\subsubsection{The Elowyn-Bonhomme map revisited}
 Let $\xi$, $F_0$ and $G_{\xi}$ be as in the original EB map (see p.~\pageref{EBmappp}) and consider the map $\mathcal{P}$ defined
 \textsl{via} (\ref{themapping2}), with
 $$\mathcal{F}_0(x,\lambda_1)=F_0(x)+\lambda_1\begin{pmatrix}-\frac{28}{5}x_1+9x_2-\frac{9}{5}x_1^2-\frac{41}{10}x_1x_2-\frac{57}{10}x_2^2 \\ -9x_1-\frac{28}{5}x_2-9x_1^2+\frac{49}{10}x_1x_2-10x_2^2 \end{pmatrix},$$
 $$\mathcal{G}_{\xi}(x,\lambda_2)=G_{\xi}(x)+\lambda_2\begin{pmatrix} 28-26x_1+\frac{2}{5}x_2+\frac{19}{2}x_1^2-\frac{33}{10}x_1x_2-\frac{2}{5}x_2^2\\ \frac{2}{5}-\frac{38}{5}x_1-6x_2+\frac{28}{5}x_1^2-\frac{7}{2}x_1x_2+\frac{23}{5}x_2^2\end{pmatrix}$$ and $$\mathcal{J}(x,\lambda)=\begin{pmatrix} \frac{67}{10}+\frac{32}{5}x_1+\frac{18}{5}x_2-\frac{33}{5}x_1^2-x_1x_2-\frac{8}{5}x_2^2 \\ \frac{59}{10}+\frac{31}{5}x_1-\frac{46}{5}x_2+\frac{5}{2}x_1^2+\frac{79}{10}x_1x_2+\frac{73}{10}x_2^2\end{pmatrix}.$$
 Then, according to proposition \ref{themappingprop2},
$\mathcal{P}$ satisfies (P1)$-$(P3). By the preceding discussion, the
bifurcation diagrams of the EB map $P$ is topologically equivalent to that of
$\mathcal{P}$ when $\lambda$ is close enough to the origin. \par As an example,
let $0.01<\omega_*=\frac{12}{37\pi}$ and define $\mathcal{P}_{0.01}$ as the
restriction of $\mathcal{P}$ to the circle $\gamma_{0.01}(s)=0.01\left(
\cos(s),\sin(s) \right)^{\!\top}$ in parameter space. Compare the diagram shown
 (on the left) in \begin{figure}[t]
\begin{center}
\includegraphics[height=170pt]{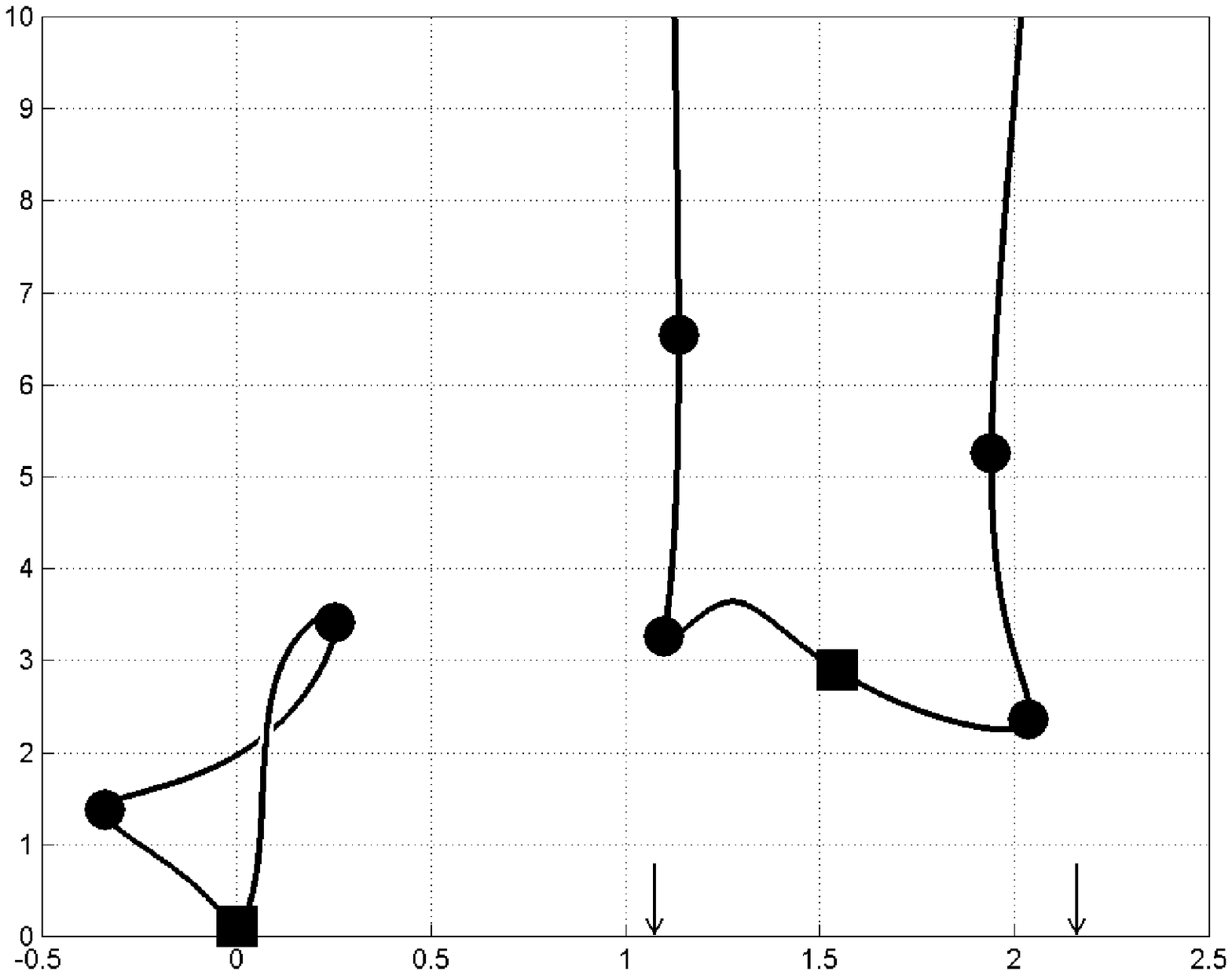}\qquad \includegraphics[height=170pt]{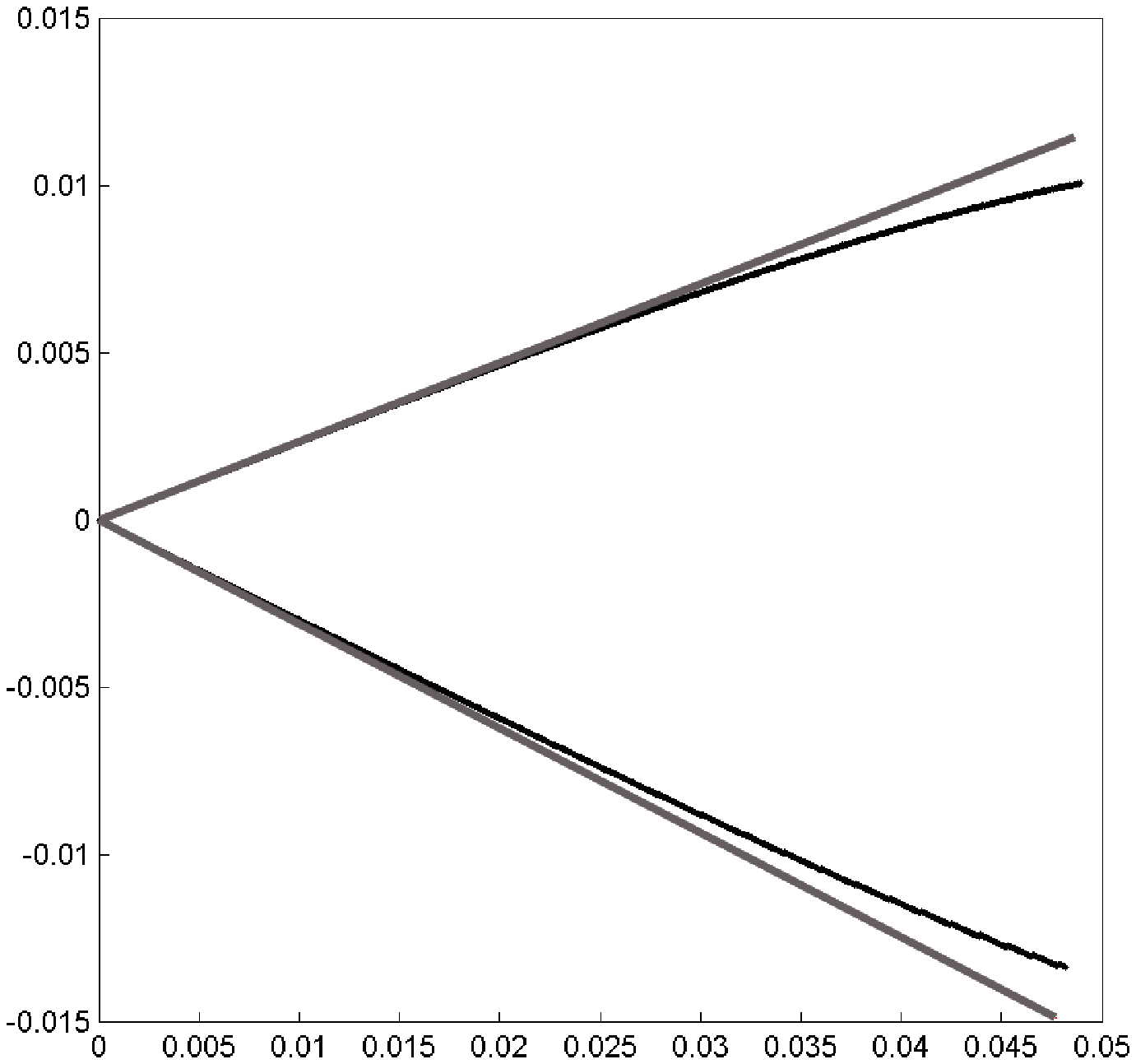}
\caption{Fixed point branches of $\mathcal{P}_{0.01}$ (left); anchoring wedges
$\mathsf{W}_{0}$ (gray) and $\mathfrak{W}_{0}$ (black) for the Elowyn-Bonhomme
maps (right).}
\end{center}\hrule \end{figure}figure~\thefigure~
with the corresponding diagram of section~\ref{EBmappp}.
\par Finally, the wedge regions
$\mathsf{W}_{0}$ of the EB map and $\mathfrak{W}_{0}$ of the revisited EB map
are shown (on the right) in
 figure~\thefigure, illustrating the last remark
of section 4.

\section*{References}
\bibliographystyle{nonlin}
\bibliography{../Bibliography/spiralsp}
\end{document}